\documentclass[12pt]{report} 
\usepackage{latexsym}
\newcommand{\insp}{\hspace*{1cm}}

\newcommand{\ds}{\displaystyle}

\usepackage{amsfonts}

\usepackage{epsfig}

\newcommand{\ad}{\mbox{ad}}

\newcommand{\reg}{\mbox{\scriptsize{reg}}}

\newcommand{\codim}{\mbox{codim}}
\newcommand{\trdeg}{\mbox{trdeg}}
\newcommand{\sing}{\mbox{\scriptsize{sing}}}

\newcommand{\rank}{\mbox{rank}}

\newcommand{\Der}{\mbox{Der}}

\newcommand{\vet}[1]{\mbox{\boldmath$#1$}}

\newcommand{\IR}{I\!\!R}

\newcommand{{\CB}}{\cal B}

\newcommand{{\CC}}{\cal C}
\def\subs{\mathop{\subset}}

 1
 1

\def\picture#1 by #2 (#3){
\vbox to #2{
\hrule width #1 height 0pt depth 0pt
\vfill
\special{picture #3}}}

\begin{document}
\begin{center}
{\bf \Large The polynomiality of the Poisson center and semi-center of a Lie algebra and Dixmier's fourth problem}\\
\ \\
{\bf \large Alfons I. Ooms}\\
\ \\
{\it Mathematics Department, Hasselt University, Agoralaan, Campus Diepenbeek, 3590 Diepenbeek, Belgium\\
E-mail address: alfons.ooms@uhasselt.be}
\end{center}

\ \\
{\bf Abstract.}\\
Let $\mathfrak{g}$ be a finite dimensional Lie algebra over an algebraically closed field $k$ of characteristic zero.  
We provide necessary and also some sufficient conditions in order for its Poisson center and semi-center to be polynomial algebras over $k$.\\
This occurs for instance  if $\mathfrak{g}$ is quadratic of index 2 with $[\mathfrak{g}, \mathfrak{g}] \neq \mathfrak{g}$ and also if $\mathfrak{g}$ is nilpotent of index at most 2.  The converse holds for filiform Lie algebras of type $L_n$, $Q_n$, $R_n$ and $W_n$.\\
We show how Dixmier's fourth problem for an algebraic Lie algebra $\mathfrak{g}$ can be reduced to that of its canonical truncation $\mathfrak{g}_\Lambda$. Moreover, Dixmier's statement holds for all Lie algebras of dimension at most eight.  
The nonsolvable, indecomposable ones among them possess a polynomial Poisson center and semi-center.\\
\ \\
{\bf \large 1. Introduction}\\
Let $\mathfrak{g}$ be a Lie algebra over an algebraically closed field $k$ of characteristic zero, with basis $x_1,\ldots, x_n$. Let $U(\mathfrak{g})$ be its enveloping algebra with center $Z(U(\mathfrak{g}))$ and semi-center $Sz(U(\mathfrak{g}))$, i.e. the subalgebra of $U(\mathfrak{g})$ 
generated by the semi-invariants of $U(\mathfrak{g})$.  Denote by $D(\mathfrak{g})$ the quotient division ring of $U(\mathfrak{g})$ with center $Z(D(\mathfrak{g}))$.  
In this paper we address the following problems:
\begin{itemize}
\item[1)] When are $Z(U(\mathfrak{g}))$ and $Sz(U(\mathfrak{g}))$ polynomial algebras over $k$ ?
\item[2)] Is $Z(D(\mathfrak{g}))$ always rational over $k$ ? (Dixmier's fourth problem [D6, p.354]).
\end{itemize}
In order to simplify things we consider the symmetric algebra $S(\mathfrak{g})$ which we identify with the polynomial algebra $k[x_1,\ldots, x_n]$.\\
We equip $S(\mathfrak{g})$ with its natural Poisson structure.  Its Poisson center $Y(\mathfrak{g})$ coincides with the algebra $S(\mathfrak{g})^{\mathfrak{g}}$ of invariants.
By a celebrated result of Michel Duflo [Du1, Du2, Du3] there exists an algebra isomorphism between $Z(U(\mathfrak{g}))$ and $Y(\mathfrak{g})$.  
Rentschler and Vergne [RV] later extended this to an algebra isomorphism between $Sz(U(\mathfrak{g}))$
and the semi-center $Sz(S(\mathfrak{g}))$, which is usually denoted by $Sy(\mathfrak{g})$.  Furthermore, $Z(D(\mathfrak{g}))$ is isomorphic 
with $R(\mathfrak{g})^{\mathfrak{g}}$, the subfield of invariants of
$R(\mathfrak{g})$, where $R(\mathfrak{g})$ is the quotient field of $S(\mathfrak{g})$.  Therefore it suffices to deal with both problems in 
$S(\mathfrak{g})$ and $R(\mathfrak{g})$, where 
things are easier and where it is possible to use MAPLE for the less trivial calculations.  Our first objective is to collect necessary (see 3.1) 
and sufficient (see 3.2) conditions in order to have polynomiality.\\
The index $i(\mathfrak{g})$ of $\mathfrak{g}$ (see 2.1) will play a  major role.  However, an alternative index $j(\mathfrak{g})$ (see 2.2) 
will perform better in the nonalgebraic case.  For instance we have the following.\\
\ \\
{\bf Theorem 1.}
$$j(\mathfrak{g}) = \mbox{trdeg}_k R(\mathfrak{g})^{\mathfrak{g}} = \mbox{trdeg}_kZ(D(\mathfrak{g})) \leq i(\mathfrak{g})$$
Moreover, equality occurs if $\mathfrak{g}$ is ad-algebraic or if $\mathfrak{g}$ has no proper semi-invariants in $S(\mathfrak{g})$.\\
\ \\
For brevity we will call $\mathfrak{g}$ coregular if $Y(\mathfrak{g})$ is a polynomial algebra over $k$.\\
\ \\
{\bf Definition 17.} Let $p_{\mathfrak{g}} \in S(\mathfrak{g})$ be the fundamental semi-invariant of $\mathfrak{g}$ (see 2.6).  We say that $\mathfrak{g}$ satisfies the Joseph-Shafrir conditions, JS
for short, if $\mathfrak{g}$ is unimodular for which $p_{\mathfrak{g}}$ is an invariant and $\mbox{trdeg}_kY(\mathfrak{g}) = i(\mathfrak{g})$. (For example JS is satisfied if $\mathfrak{g}$ has no proper
semi-invariants (Remark 2)).\\
\ \\
The following criterion for coregularity will be employed quite often.\\
\ \\
{\bf Corollary 19.} (Short version)\\
Let $\mathfrak{g}$ be a Lie algebra satisfying JS with center $Z(\mathfrak{g})$.  If $\mathfrak{g}$ is coregular then
$$3i(\mathfrak{g}) + 2\ \mbox{deg}\ p_{\mathfrak{g}} \leq \mbox{dim}\ \mathfrak{g} + 2\ \mbox{dim}\ Z(\mathfrak{g})$$
This inequality imposes a strong upperbound on $i(\mathfrak{g})$.  Therefore coregularity becomes a rare phenomenon for nonabelian Lie algebras having a large index.  
This is especially true if $i(\mathfrak{g}) = \mbox{dim}\ \mathfrak{g} - 2$ as it is in the following result (see also Proposition 40 and Corollary 41).\\
\ \\
{\bf Theorem 20.} Let $\mathfrak{g}$ be unimodular, having an abelian ideal $\mathfrak{h}$ of codimension one.  Then the following are equivalent:
\begin{itemize}
\item[(1)] $\mathfrak{g}$ is coregular and $\mbox{trdeg}_kY(\mathfrak{g}) = i(\mathfrak{g})$
\item[(2)] $\mathfrak{g}$ is coregular and ad-algebraic
\item[(3)] $\dim[\mathfrak{g},\mathfrak{g}] \leq 2$
\end{itemize}
The conditions in the above results cannot be weakened as shown by the examples 24-26.\\
As an application of this theorem we consider the nilradical $N$ of the parabolic subalgebra $P$ of type $(1, 1, n-2)$ inside $sl(n)$ and show that
\begin{center}
$N$ is not coregular $\ \ \ \Leftrightarrow \ \ \ n \geq 5$
\end{center}
(see Proposition 28).  This extends an example by A. Hersant [J2, 8.5].\\
\ \\
At the end of section 3 we consider a Lie algebra $\mathfrak{g}$ for which $Y(\mathfrak{g})$ is saturated with quotient field $R(\mathfrak{g})^{\mathfrak{g}}$.  
Then there exist irreducible, proper semi-invariants $v_1, \ldots, v_t \in S(\mathfrak{g})$ such  that $Sy(\mathfrak{g}) = Y(\mathfrak{g})[v_1,\ldots, v_t]$ 
is a polynomial ring over $Y(\mathfrak{g})$. In particular, if $Y(\mathfrak{g})$ is polynomial over $k$, then so is $Sy(\mathfrak{g})$ (Theorem 36).\\
\ \\
The above conditions are satisfied if $j(\mathfrak{g}) = \dim Z(\mathfrak{g})$ because then $Y(\mathfrak{g}) = S(Z(\mathfrak{g}))$ and 
$R(\mathfrak{g})^{\mathfrak{g}} = R(Z(\mathfrak{g}))$.  
Hence both $Y(\mathfrak{g})$ and $Sy(\mathfrak{g})$ are polynomial over $k$, while $R(\mathfrak{g})^{\mathfrak{g}}$ is rational over $k$.  
Moreover, $v_1, \ldots, v_t$ are then precisely the irreducible
factors of a special semi-invariant $p'_\mathfrak{g} \in S(\mathfrak{g})$, which takes over the role of the fundamental semi-invariant $p_\mathfrak{g}$ (Theorem 37).  
This is illustrated in Example 39 and applied to some Lie algebras such as $L_{8,25}$
(see Example 58) of section 5.\\
\ \\
In section 4 we study the coregularity for Lie algebras of index at most two.  The following is one of the main results:\\
\ \\
{\bf Theorem 45.} Any nilpotent Lie algebra with index at most two is coregular.\\
\ \\
Its proof is constructive and we can give a useful characterization of the generator(s) of $Y(\mathfrak{g})$. (see the claim within the proof).  
This is used in the following application to the major types of filiform Lie algebras:\\
\ \\
{\bf Theorem 51.}
\begin{itemize}
\item[(1)] If $\mathfrak{g}$ is of type $Q_n$ or $W_n$ then $\mathfrak{g}$ is coregular since $i(\mathfrak{g}) \leq 2$.
\item[(2)] If $\mathfrak{g}$ is of type $L_n$ or $R_n$ then
\end{itemize}
\begin{center}
$\mathfrak{g}$ is coregular $\ \ \ \Leftrightarrow \ \ \ i(\mathfrak{g}) \leq 2$
\end{center}

Theorem 45 cannot be extended to the solvable case as there exists a solvable Lie algebra of index two which is not coregular (Example 23.).\\
\ \\
{\bf Theorem 52.} Let $\mathfrak{g}$ be a quadratic Lie algebra.  Then $\mathfrak{g}$ is coregular if one of the following conditions is satisfied:
\begin{itemize}
\item[(i)] $[\mathfrak{g},\mathfrak{g}]\neq \mathfrak{g}$ and $i(\mathfrak{g}) = 2$
\item[(ii)] $\mathfrak{g}$ is nilpotent and $i(\mathfrak{g}) = 3$
\end{itemize}
\ \\
In section 5 we verify, case by case, that any nonsolvable, indecomposable Lie algebra of dimension at 
most eight satisfies the two problems we raised in the beginning (Theorem 53).\\
However, in dimension nine, we exhibit a counterexample (Example 59).\\
\ \\
Section 6 is devoted to Dixmier's fourth problem.\\
We list some important classes where this question is known to have a positive answer and we prove that it is also the case for all Lie algebras of dimension at most 8 
(Proposition 63).  The following is the main result of this section:\\
\ \\
{\bf Theorem 66.} Let $\mathfrak{g}$ be an algebraic Lie algebra for which the field $Z(D(\mathfrak{g}_\Lambda))$ is 
freely generated by semi-invariants $u_1,\ldots, u_s$ of $U(\mathfrak{g})$.  Then $Z(D(\mathfrak{g}))$ is rational over $k$.\\
\ \\
As an application we obtain a result by Panyushev [Pa1], namely $Z(D(\mathfrak{g}))$ is rational over $k$ if $\mathfrak{g}$ is 
any biparabolic subalgebra of a simple Lie algebra of type $A$ or $C$ (Corollary 67).\\
\ \\
Some of the results of [AOV2] are used in section 5.  Therefore we briefly discuss the well known Gelfand-Kirillov conjecture.  
This is a much stronger statement than Dixmier's
fourth problem (see also Example 60 and Proposition 62).  In the Appendix we correct the proof of an example from [GK] showing the existence of 
nonalgebraic Lie algebras satisfying the Gelfand-Kirillov conjecture (Example 71).\\
\ \\
{\bf \large 2. Preliminaries}\\
{\bf 2.1 $i(\mathfrak{g})$, the index of $\mathfrak{g}$}\\
Let $k$ be an algebraically closed field of characteristic zero and let $\mathfrak{g}$ be a Lie algebra over $k$ with basis $x_1, \ldots, x_n$. For each $\xi \in \mathfrak{g}^\ast$ we consider its stabilizer
\begin{eqnarray*}
\mathfrak{g}(\xi) = \{x \in \mathfrak{g}\mid \xi ([x,y]) = 0\ \mbox{for all}\ y \in \mathfrak{g}\}
\end{eqnarray*}
The minimal value of $\dim \mathfrak{g}(\xi)$ is called the index of $\mathfrak{g}$ and is denoted by $i(\mathfrak{g})$ [D6, 1.11.6; TY, 19.7.3].
Put $c(\mathfrak{g}) = (\dim \mathfrak{g} + i(\mathfrak{g}))/2$.  This integer will play an important role throughout this paper.  
An element $\xi \in \mathfrak{g}^\ast$ is called regular if $\dim \mathfrak{g}(\xi) = i(\mathfrak{g})$.  
The set $\mathfrak{g}^\ast_{\reg}$ of all regular elements of $\mathfrak{g}^\ast$ is an open dense subset of $\mathfrak{g}^\ast$.\\
We put $\mathfrak{g}^\ast_{\sing} = \mathfrak{g}^\ast \backslash \mathfrak{g}^\ast_{\reg}$.  
Clearly, $\codim\ \mathfrak{g}^\ast_{\sing} \geq 1$.  Following [JS] we call $\mathfrak{g}$ singular if equality holds and nonsingular otherwise.  
For instance, any semi-simple Lie algebra $\mathfrak{g}$ is nonsingular since $\codim\ \mathfrak{g}^\ast_{\sing} = 3$. We recall from [D6, 1.14.13] that
$$i(\mathfrak{g}) = \dim \mathfrak{g} - \rank_{R(\mathfrak{g})} ([x_i, x_j])$$
In particular, $\dim \mathfrak{g} - i(\mathfrak{g})$ is an even number.\\
\ \\
{\bf 2.2 $j(\mathfrak{g})$, the alternative index of $\mathfrak{g}$}\\
Let $H$ be the algebraic hull of $\ad~\mathfrak{g}$ in $\mbox{Der}~\mathfrak{g}$ ([C, p.173; TY, 24.5.4]), i.e. the smallest algebraic Lie subalgebra $H$ of 
$\mbox{Der}\ \mathfrak{g}$ containing $\ad~\mathfrak{g}$.
Let $\xi \in \mathfrak{g}^\ast$ and put
$$\mathfrak{g}[\xi] = \{x \in \mathfrak{g}\mid \xi (Ex)) = 0\ \mbox{for all}\ E \in H\}$$
This is an ideal of $\mathfrak{g}(\xi)$ which contains the center $Z(\mathfrak{g})$ of $\mathfrak{g}$.  Clearly, $\mathfrak{g}[\xi] = \mathfrak{g}(\xi)$ if 
$\mathfrak{g}$ is ad-algebraic (i.e. $\ad~\mathfrak{g} = H$).  Let $E_1,\ldots, E_m$ be a basis of $H$.  Then it is easily seen that
$$\dim \mathfrak{g}[\xi] = \dim \mathfrak{g} - \mbox{rank}(\xi(E_i x_j))$$
We denote by $j(\mathfrak{g})$ the minimal value of $\dim \mathfrak{g}[\xi]$, $\xi \in \mathfrak{g}^\ast$.  Then
$$j(\mathfrak{g}) = \dim \mathfrak{g} - \mbox{rank}_{R(\mathfrak{g})}(E_ix_j)$$
Clearly, $\dim Z(\mathfrak{g}) \leq j(\mathfrak{g}) \leq i(\mathfrak{g})$ and $j(\mathfrak{g}) = i(\mathfrak{g})$ if $\mathfrak{g}$ is ad-algebraic.  For the first part of the following we refer to [O1, O2; 
RV, p.401].\\
\ \\
{\bf Theorem 1.}
$$j(\mathfrak{g}) = \mbox{trdeg}_kR(\mathfrak{g})^{\mathfrak{g}}= \mbox{trdeg}_kZ(D(\mathfrak{g})) \leq i(\mathfrak{g})$$
Moreover, equality occurs if one of the following conditions is satisfied:
\begin{itemize}
\item[(1)] $\mathfrak{g}$ is ad-algebraic
\item[(2)] $\mathfrak{g}$ has no proper semi-invariants in $S(\mathfrak{g})$ (or equivalently in $U(\mathfrak{g})$) [OV, Proposition 4.1].
\end{itemize}

{\bf 2.3 Commutative polarizations of $\mathfrak{g}$}\\
Suppose $\mathfrak{g}$ admits a commutative Lie subalgebra $\mathfrak{h}$ such that $\dim \mathfrak{h} = c(\mathfrak{g})$, i.e. $\mathfrak{h}$ is a commutative polarization 
(notation: CP) with respect to any $\xi \in \mathfrak{g}_{\reg}^\ast$ [D6, 1.12].\\
These CP's occur frequently in the nilpotent case [O7, O8].  If in addition $\mathfrak{h}$ is an ideal of $\mathfrak{g}$ then we call $\mathfrak{h}$ a CP-ideal (notation: CPI).
If a solvable Lie algebra $\mathfrak{g}$ admits a CP then it also admits a CPI [EO, Theorem 4.1]. \\
\ \\
{\bf 2.4 The Poisson algebra $S(\mathfrak{g})$ and its center}\\
The symmetric algebra $S(\mathfrak{g})$, which we identify with $k[x_1,\ldots, x_n]$, has a natural Poisson algebra structure, the Poisson bracket of $f, g \in S(\mathfrak{g})$ given by:
$$\{f,g\} = \sum\limits_{i=1}^n \sum\limits_{j=1}^n [x_i,x_j] \ds\frac{\partial f}{\partial x_i} \ds\frac{\partial g}{\partial x_j}$$
In particular, $S(\mathfrak{g}),\{,\}$ is a Lie algebra for which $\mathfrak{g}$ is a Lie subalgebra since for any two elements $x, y \in \mathfrak{g}$ we have that $\{x,y\} = [x,y]$.  Also, for all $f,g,h \in S(\mathfrak{g})$:
$$\{f, gh\} = \{f,g\}h + g\{f,h\} \insp (\ast)$$
It now easily follows that the center of $S(\mathfrak{g}),\{,\}$ is equal to
$$\{f \in S(\mathfrak{g})\mid \{x,f\} = 0\ \ \forall x \in \mathfrak{g}\}$$
and since $\{x,f\} = \ad~x(f)$ this clearly coincides with $Y(\mathfrak{g}) = S(\mathfrak{g})^{\mathfrak{g}}$, the subalgebra of invariant polynomials of 
$S(\mathfrak{g})$.\\
The Poisson bracket has a unique extension to the quotient field $R(\mathfrak{g})$ of $S(\mathfrak{g})$ such that $(\ast)$ holds in $R(\mathfrak{g})$.  
It follows that $R(\mathfrak{g}),\{,\}$ is a Lie algebra with center $R(\mathfrak{g})^{\mathfrak{g}}$, the subfield of rational invariants of 
$R(\mathfrak{g})$. $R(\mathfrak{g})$ is called the rational Poisson algebra [V, p. 311].\\
\ \\
Let $A$ be a Poisson commutative subalgebra of $S(\mathfrak{g})$ (i.e. $\{f,g\} = 0$ for all $f,g \in A$).  Then it is well-known that $\trdeg_k(A) \leq c(\mathfrak{g})$.  
$A$ is called complete if equality holds and strongly complete if it is also a maximal Poisson commutative subalgebra.  According to Sadetov there always exists a complete 
Poisson commutative subalgebra of $S(\mathfrak{g})$ [Sa].  For example, suppose $\mathfrak{g}$ admits a commutative polarization (CP) $\mathfrak{h}$.  
Then $S(\mathfrak{h})$ is a polynomial, strongly complete subalgebra of $S(\mathfrak{g})$ and its quotient field $R(\mathfrak{h})$ is a maximal Poisson commutative 
subfield of $R(\mathfrak{g})$ [O4, Theorem 14].\\
\ \\
{\bf 2.5 The semi-center $Sy(\mathfrak{g})$ of $S(\mathfrak{g})$}\\
Let $\lambda \in \mathfrak{g}^\ast$.  We denote by $S(\mathfrak{g})_{\lambda}$ the set of all $f \in S(\mathfrak{g})$ such that $\ad~x(f) = \lambda(x)f$ for all 
$x \in \mathfrak{g}$.  Any element $f \in S(\mathfrak{g})_{\lambda}$ is said to be a semi-invariant w.r.t. the weight $\lambda$.  We call $f$ a proper semi-invariant 
if $\lambda \neq 0$.  Clearly, $S(\mathfrak{g})_{\lambda} S(\mathfrak{g})_{\mu} \subset S(\mathfrak{g})_{\lambda + \mu}$ for all $\lambda, \mu \in \mathfrak{g}^\ast$.  
Let $f, g \in S(\mathfrak{g})$.  If $fg$ is a nonzero semi-invariant of $S(\mathfrak{g})$, then so are $f$ and $g$.\\
The sum of all $S(\mathfrak{g})_{\lambda}$, $\lambda \in \mathfrak{g}^\ast$, is direct and it is a nontrivial factorial subalgebra $Sy(\mathfrak{g})$ of $S(\mathfrak{g})$ [D3, Mo, LO]. Moreover, it is Poisson commutative [OV, p. 308].\\
Any nonzero semi-invariant can be written uniquely as a product of irreducible semi-invariants.\\
Suppose $h \in R(\mathfrak{g})$, $h \neq 0$.  Then $h \in R(\mathfrak{g})^{\mathfrak{g}}$ if and only if $h$ can be written as a quotient of two semi-invariants of the same weight.\\
\ \\
{\bf Remark 2.} Assume that $\mathfrak{g}$ has no proper semi-invariants (as it is if the radical of $\mathfrak{g}$ is nilpotent).  Then $R(\mathfrak{g})^{\mathfrak{g}}$ is the quotient field of $S(\mathfrak{g})^{\mathfrak{g}} = Y(\mathfrak{g})$.  In particular,
$$\trdeg_k Y(\mathfrak{g}) = \trdeg_k R(\mathfrak{g})^{\mathfrak{g}} = i(\mathfrak{g})$$
by Theorem 1.  Also, $\mathfrak{g}$ is unimodular (i.e. $\mbox{tr}(\ad~x) = 0$ for all $x \in \mathfrak{g}$) by [DDV, Thm. 1.11] and its proof.\\
\ \\
The weights of the semi-invariants of $S(\mathfrak{g})$ form an additive semi-group $\Lambda(\mathfrak{g})$, which is not necessarily finitely generated [DDV, p. 322].  
However, the subgroup $\Lambda_R(\mathfrak{g})$  of $\mathfrak{g}^\ast$ generated by $\Lambda(\mathfrak{g})$ is a finitely generated free abelian group [NO, Theorem 1.3], [FJ2, p. 1519].\\
Next, we denote by $\mathfrak{g}_\Lambda$ the intersection of $\mbox{ker} \lambda$, $\lambda \in \Lambda(\mathfrak{g})$.  $\mathfrak{g}_\Lambda$ is a characteristic ideal of 
$\mathfrak{g}$ which contains $[\mathfrak{g},\mathfrak{g}]$.  It is called the canonical truncation of $\mathfrak{g}$.\\
\ \\
{\bf Lemma 3.} Let $u \in S(\mathfrak{g})$ be a nonzero semi-invariant with weight $\lambda \in \Lambda (\mathfrak{g})$.  Denote by $C(u)$ (resp. $C_R(u)$) 
the centralizer of $u$ in
$S(\mathfrak{g})$ (resp. $R(\mathfrak{g})$).  Then we have
$$C(u) = S(\mbox{ker} \lambda)\ \ \ \mbox{and}\ \ \ C_R(u) = R(\mbox{ker} \lambda)$$

{\bf Proof.} We may assume that $\lambda \neq 0$.  Choose a basis $x_1, x_2,\ldots, x_n$ of $\mathfrak{g}$ such that $x_2,\ldots, x_n$ is a basis of $\mbox{ker}\lambda$ and such 
that $\lambda(x_1) = 1$.  Since $u \in S(\mathfrak{g})_\lambda$ we have $\{f,u\} = d_\lambda(f)u$ for all $f \in S(\mathfrak{g})$ [OV, p.308].  Then the first equality follows from
$$f \in C(u) \ \ \ \Leftrightarrow \ \ \ d_\lambda (f) = 0 \ \ \ \Leftrightarrow \ \ \ \ds\frac{\partial f}{\partial x_1} = 0 \ \ \ \Leftrightarrow \ \ \ f \in k [x_2,\ldots, x_n] = S(\mbox{ker} \lambda)$$
Next, we take a nonzero $h \in C_R(u)$.  We may write $h = f/g$ for some nonzero, relatively prime $f,g \in S(\mathfrak{g})$.  From $hg = f$ we deduce $h \{g,u\} = \{f,u\}$ 
since $\{h,u\} = 0$. Hence,
$$hd_\lambda(g)u = d_\lambda(f)u$$
Simplification gives
$$d_\lambda(g) f = d_\lambda(f)g$$
Now suppose $d_\lambda(f) \neq 0$.  Then $f$, being coprime with $g$, divides $d_\lambda(f)$, contradicting the fact that $\mbox{deg}\ d_\lambda(f) < \mbox{deg}\ f$.\\
Therefore $d_\lambda(f) = 0$ and thus $f \in S(\mbox{ker} \lambda)$. Similarly, $g \in S(\mbox{ker} \lambda)$ and so $h = f/g \in R(\mbox{ker} \lambda)$.  
Consequently, $C_R(u) \subset R(\mbox{ker} \lambda)$. The other inclusion is obvious. \hfill $\square$\\
\ \\
Using this lemma one can now apply the same approach as in [DNO, pp. 331-334] and [MO, pp. 213-214] in order to obtain the following.   
In fact (1), (2), (3) do not require for $k$ to be algebraically closed.  See also [BGR, F, FJ2, RV].\\
\ \\
{\bf Theorem 4.}
\begin{itemize}
\item[1.] $C(Sy(\mathfrak{g})) = S(\mathfrak{g}_\Lambda)$ and $ C_R(Sy(\mathfrak{g})) = R(\mathfrak{g}_\Lambda)$
\item[2.] $\mathfrak{g}_\Lambda$ has no proper semi-invariants and so $R(\mathfrak{g}_\Lambda)^{\mathfrak{g}_\Lambda}$ is the quotient field of $Y(\mathfrak{g}_\Lambda)$.
Also $\mbox{trdeg}_kY(\mathfrak{g}_\Lambda) = i(\mathfrak{g}_\Lambda)$
\item[3.] $S(\mathfrak{g})^{\mathfrak{g}_\Lambda} = Y(\mathfrak{g}_\Lambda)$ and $R(\mathfrak{g})^{\mathfrak{g}_\Lambda} = R(\mathfrak{g}_\Lambda)^{\mathfrak{g}_\Lambda}$
\item[4.] $c(\mathfrak{g}_\Lambda) = c(\mathfrak{g})$ (use [OV, Lemma 3.7] and [O7, Proposition 3.2])
\item[5.] $Sy(\mathfrak{g}) \subset Y(\mathfrak{g}_\Lambda) = Sy(\mathfrak{g}_\Lambda)$ and equality occurs if $\mathfrak{g}$ is almost algebraic or if $\mathfrak{g}$ 
is Frobenius (i.e. $i(\mathfrak{g}) = 0$)
\item[6.] Suppose $\mathfrak{h}$ is a CP-ideal of $\mathfrak{g}$. Then $\mathfrak{h} \subset \mathfrak{g}_\Lambda$ [EO, p. 141] and $Y(\mathfrak{g}_\Lambda) 
\subset S(\mathfrak{h})$
\end{itemize}

{\bf 2.6 The fundamental semi-invariant $p_\mathfrak{g}$}\\
{\bf Definition 5.} Put $t = \dim \mathfrak{g} - i(\mathfrak{g})$, which is the rank of the structure matrix $B = ([x_i,x_j]) \in M_n(R(\mathfrak{g}))$, 
where $x_1, \ldots, x_n$ is an arbitrary basis of $\mathfrak{g}$.  Assume first that $\mathfrak{g}$ is nonabelian.  Then the greatest common divisor 
$q_{\mathfrak{g}}$ of the $t \times t$ minors in $B$ is a nonzero semi-invariant of $S(\mathfrak{g})$ [DNO, pp. 336-337].  If $\mathfrak{g}$ is abelian we put 
$q_{\mathfrak{g}} = 1$.  Next, let $p_{\mathfrak{g}}$ be the greatest common divisor of the Pfaffians of the principal $t \times t$ minors in $B$.  
In particular, $\deg p_{\mathfrak{g}} \leq (\dim \mathfrak{g} - i(\mathfrak{g}))/2$.  By [OV, Lemma 2.1] $p_{\mathfrak{g}}^2 = q_{\mathfrak{g}}$ up to a 
nonzero scalar multiplier.  We call $p_{\mathfrak{g}}$ the fundamental semi-invariant of $S(\mathfrak{g})$ (instead of $q_{\mathfrak{g}}$ as we did in [OV, p. 309]).\\
\ \\
{\bf Remark 6.} [OV, p. 307]
\begin{center}
$\mathfrak{g}$ is singular if and only if $p_{\mathfrak{g}} \notin k$
\end{center}

{\bf Example 7.}  Let $\mathfrak{g}$ be a nonabelian Lie algebra with center $Z(\mathfrak{g})$.  $\mathfrak{g}$ is called square integrable (SQ.I.) if 
$i(\mathfrak{g}) = \dim Z(\mathfrak{g})$. For instance any Heisenberg Lie algebra is square integrable.\\
Choose a basis $x_1, \ldots, x_t, x_{t+1},\ldots, x_n$ such that $x_{t+1},\ldots, x_n$ is a basis of $Z(\mathfrak{g})$.\\
Then, $t = \dim \mathfrak{g} - \dim Z(\mathfrak{g}) = \dim\mathfrak{g} - i(\mathfrak{g})$, which is the rank of the matrix $([x_i, x_j])_{1\leq i,j\leq t}$.  
By the above, its Pfaffian coincides with $p_{\mathfrak{g}}$ (up to a nonzero scalar).  Hence, $\deg p_{\mathfrak{g}} = (\dim\mathfrak{g} - i(\mathfrak{g}))/2 \geq 1$ 
and so $\mathfrak{g}$ is singular. In particular, any Frobenius Lie algebra $\mathfrak{g}$ is singular.\\
\ \\
{\bf Lemma 8.} [J5, Lemma 2.3] Let $\mathfrak{g}$ be an algebraic Lie algebra. Then $p_{\mathfrak{g}_\Lambda}$ divides $p_{\mathfrak{g}}$.\\
\ \\
{\bf 2.7 Frobenius Lie algebras}\\
A Lie algebra $\mathfrak{g}$ is called Frobenius if there is a linear functional $\xi \in \mathfrak{g}^\ast$ such that the alternating bilinear 
$B_\xi (x,y) = \xi([x,y]), x, y \in \mathfrak{g}$, is nondegenerate, i.e. $i(\mathfrak{g}) = 0$.  The name was suggested to us by George Seligman because of its obvious resemblance
with the notion of an associative Frobenius algebra.  They came about in connection with Jacobson's problem on the characterization of Lie algebras having a primitive universal enveloping algebra.
It turns out that:\\
$U(\mathfrak{g})$ is primitive if and only if $Z(D(\mathfrak{g})) = k$, i.e. $j(\mathfrak{g}) = 0$ [O1, O2].\\
In particular, $U(\mathfrak{g})$ is  primitive if $\mathfrak{g}$ is Frobenius and the converse holds if $\mathfrak{g}$ is ad-algebraic by Theorem 1.\\
Frobenius Lie algebras form a large class and they appear naturally in different areas.  For example many parabolic and biparabolic (seaweed) subalgebras of semi-simple Lie algebras are Frobenius 
[CGM, CMW, DY, E1, E2, E3, CV, JS, PY2, O3], including most Borel subalgebras of simple Lie algebras [EO, p. 146].  
A Frobenius biparabolic Lie algebra $\mathfrak{g}$ satisfies interesting properties.  
For instance $\mathfrak{g}_\Lambda = [\mathfrak{g},\mathfrak{g}]$ [J5, Proposition 7.6], which is not true for all Frobenius Lie algebras as the following demonstrates 
(this answers a question by Joseph [J5, Remark 7.6].\\
\ \\
{\bf Example 9.} Let $L$ be the Lie algebra over $k$ with basis $x_1, x_2, x_3, x_4$ and nonvanishing brackets $[x_1, x_3] = x_3, \ [x_1, x_4] = x_4,\ [x_2, x_3] = x_4$.\\
Consider its structure matrix $B = ([x_i, x_j])$.  Clearly $\mbox{det} B = x_4^4 \neq 0$.  Hence $i(L) = 0$ by 2.1 and $p_L = x_4^2$.\\
$x_4$ is the only irreducible semi-invariant of $S(L)$ (see below).  Its weight $\lambda \in L^\ast$ is determined by $\lambda(x_1) = 1, \lambda (x_2) = \lambda (x_3) = \lambda(x_4) = 0$.\\
Consequently, $L_\Lambda = \mbox{ker} \lambda = \langle x_2, x_3, x_4\rangle$, while $[L,L] = \langle x_3, x_4\rangle$ (which happens to be  a CPI of $L$).  Moreover $Sy(L) = k[x_4] = Y(L_\Lambda)$.\\
\ \\
We now collect some useful facts on semi-invariants from [O3, DNO].  Let $\mathfrak{g}$ be a Frobenius Lie algebra with basis $x_1,\ldots, x_n$.  
Then $n$ is even and $\mathfrak{g}$ has a trivial center.  The Pfaffian $Pf([x_i, x_j]) \in S(\mathfrak{g})$ is homogeneous of degree $\ds\frac{1}{2} \mbox{dim}\ \mathfrak{g}$ 
and $(Pf([x_i, x_j]))^2 = \mbox{det} ([x_i, x_j]) \neq 0$ by 2.1.  Hence $p_\mathfrak{g} = Pf([x_i, x_j])$.  We put $\Delta(\mathfrak{g}) = \mbox{det} ([x_i,x_j])$ (which is well determined up to
nonzero scalar multipliers).\\
$p_\mathfrak{g}$ is a semi-invariant with weight $\tau$, where $\tau(x) = tr(\ad~x)$, $x \in \mathfrak{g}$.\\
Moreover, any semi-invariant of $S(\mathfrak{g})$ is homogeneous.  It is also a semi-invariant under the action of $\mbox{Der}\ \mathfrak{g}$.\\
\ \\
{\bf Theorem 10.} Let $\mathfrak{g}$ be Frobenius.  Decompose $p_\mathfrak{g}$ into a product of irreducible factors:
$$p_\mathfrak{g} = v_1^{m_1} \ldots, v_r^{m_r}, \ \ \ m_i \geq 1$$
Then:
\begin{itemize}
\item[(1)] $v_1,\ldots,v_n$ are the only (up to nonzero scalars) irreducible semi-invariants of $S(\mathfrak{g})$, say with weights $\lambda_1,\ldots, \lambda_r \in \Lambda(\mathfrak{g})$.
\item[(2)] $Sy(\mathfrak{g}) = k[v_1, \ldots, v_r] = Y(\mathfrak{g}_\Lambda)$, a polynomial algebra over $k$.
\item[(3)] $r = \dim\mathfrak{g} - \dim\mathfrak{g}_\Lambda = i(\mathfrak{g}_\Lambda)$
\item[(4)] $\lambda_1, \ldots, \lambda_r$ are linearly independent over $k$. They generate the semi-group $\Lambda(\mathfrak{g})$ and $\mathfrak{g}_\Lambda = \cap 
\ \mbox{ker} \lambda_i$, $i = 1,\ldots, r$
\item[(5)] $\mathfrak{h}_i = \mbox{ker} \lambda_i$ is an ideal of $\mathfrak{g}$ of index one and $Y(\mathfrak{h}_i) = k[v_i]$ and $R(\mathfrak{h}_i)^{\mathfrak{h}_i} = k(v_i)$
\item[(6)] $m_1 \lambda_1 + \ldots + m_r \lambda_r = \tau (\ast)$ and $m_1 \deg v_1 + \ldots + m_r \deg v_r = \deg p_\mathfrak{g} = \ds\frac{1}{2} \dim \mathfrak{g} (\ast\ast)$
\item[(7)] (Joseph [J5, 2.2]) Suppose in addition that $\mathfrak{g}$ is algebraic.  \\
Then $p_{\mathfrak{g}_\Lambda} = v_1^{m_1-1} \ldots v_r^{m_r-1}$. In particular,\\
$\mathfrak{g}_\Lambda$ is nonsingular $\ \ \ \Leftrightarrow \ \ \ m_i = 1$ for all $i = 1,\ldots, r$.
\end{itemize}

{\bf Remark 11.} Each semi-invariant $v_i$ is determined by its weight $\lambda_i$ (up to a nonzero scalar multiplier) 
[Indeed, suppose that also $v \in S(\mathfrak{g})_{\lambda_i}$, $v\neq 0$.  Then $vv_i^{-1} \in R(\mathfrak{g})^{\mathfrak{g}}= k$, i.e. $v=av_i$ for some nonzero $a \in k$].
  Therefore $\lambda_i$ will provide information on $v_i$.  For example its multiplicity $m_i$ (by $(\ast)$ since $\lambda_1,\ldots, \lambda_r$ are linearly independent over $k$), 
$\deg v_i$ and $(\ast\ast)$ can be obtained directly from $(\ast)$.  To demonstrate this we take $\xi \in \mathfrak{g}_{reg}^\ast$, i.e. $\xi (p_\mathfrak{g}) \neq 0$ (we extend $\xi$ to an algebra
endomorphism of $S(\mathfrak{g})$) and hence  also $\xi(v_i) \neq 0$.  By [O3, p. 21] there exists a unique element $x_\xi \in \mathfrak{g}$ such that $\xi \circ \ad~x_\xi = \xi$ (Nowadays $x_\xi$ is called a 
principal element of $\mathfrak{g}$). From $\ad~x_\xi(v_i) = \lambda_i(x_\xi) v_i$ we get $\xi(\ad~x_\xi (v_i)) = \lambda_i (x_\xi) \xi(v_i)$, which we can rewrite as 
$(\deg v_i) \xi (v_i) = \lambda_i(x_\xi)\xi(v_i)$ since $v_i$ is homogeneous.
Simplification yields $\deg v_i = \lambda_i(x_\xi)$.  On the other hand, $\tau(x_\xi) = tr(\ad~x_\xi) = \ds\frac{1}{2} \dim\mathfrak{g}$ [O3, Theorem 3.3]. 
Substitution in $(\ast)$ gives us $(\ast\ast)$. \hfill $\square$\\
\ \\
{\bf 2.8 The Frobenius semi-radical $F(\mathfrak{g})$}\\
Put $F(\mathfrak{g}) = \sum\limits_{\xi \in \mathfrak{g}_{\reg}^{\ast}} \mathfrak{g}(\xi)$.  
This is a characteristic ideal of $\mathfrak{g}$ containing $Z(\mathfrak{g})$ and for which $F(F(\mathfrak{g})) = F(\mathfrak{g})$.  
It can also be characterized as follows: $R(\mathfrak{g})^{\mathfrak{g}} \subset R(F(\mathfrak{g}))$ and if $\mathfrak{g}$ is algebraic then $F(\mathfrak{g})$ is 
the smallest Lie subalgebra of $\mathfrak{g}$ with this property.  Similar results hold in $D(\mathfrak{g})$  [O5 Proposition 2.4, Theorem 2.5] 
Also, $F(\mathfrak{g}) \subset \mathfrak{g}_\Lambda$.\\
As a special case we have the following:\\
\ \\
{\bf Remark 12.} $Y(\mathfrak{g}) \subset S(F(\mathfrak{g}))$ (respectively $Z(U(\mathfrak{g})) \subset U(F(\mathfrak{g})))$ and $F(\mathfrak{g})$ is the smallest Lie 
subalgebra of $\mathfrak{g}$ with this property in case $\mathfrak{g}$ is an algebraic Lie algebra without proper semi-invariants.\\
\ \\
In case $\mathfrak{g}$ is square integrable we notice that $F(\mathfrak{g}) = Z(\mathfrak{g})$ (since $\mathfrak{g}(\xi) = Z(\mathfrak{g})$ for all 
regular $\xi \in \mathfrak{g}^\ast$) which forces $R(\mathfrak{g})^{\mathfrak{g}} = R(Z(\mathfrak{g}))$.  See also Remark 38.  In particular, $Y(\mathfrak{g}) = S(Z(\mathfrak{g}))$, which is a polynomial algebra.\\
\ \\
If $\mathfrak{g}$ admits a CP $\mathfrak{h}$ then $F(\mathfrak{g})$ is commutative (since $F(\mathfrak{g}) \subset \mathfrak{h})$.  Clearly,
\begin{center}
$F(\mathfrak{g}) = 0$ if and only if $\mathfrak{g}$ is Frobenius\\
\end{center}
For this reason $F(\mathfrak{g})$ is called the Frobenius semi-radical of $\mathfrak{g}$.  
At the other end of the spectrum we have the Lie algebras for which $F(\mathfrak{g}) = \mathfrak{g}$, which we call quasi-quadratic.  
These are unimodular and they do not possess any proper semi-invariants.  They form a large class, which include all quadratic Lie algebras (and hence all 
abelian and semi-simple Lie algebras) [O5].\\
\ \\
{\bf 3. General results}\\
{\bf 3.1 Necessary conditions for polynomiality}\\
\ \\
{\bf Theorem 13.} [OV, Theorem 1.1] Let $\mathfrak{g}$ be a Lie algebra for which the semi-center $Sy(\mathfrak{g})$ is freely generated by homogeneous elements 
$f_1, \ldots, f_r$.\\
Then
$$\sum\limits_{i=1}^r \deg f_i \leq c(\mathfrak{g})$$

{\bf Definition 14.} A Lie algebra $\mathfrak{g}$ is called coregular if $Y(\mathfrak{g})$ is a polynomial algebra over $k$.\\
\ \\
{\bf Proposition 15.} [OV, Proposition 1.6].  Assume that $\mathfrak{g}$ is nonabelian, without proper semi-invariants.  
If $\mathfrak{g}$ is coregular then $\codim \ \mathfrak{g}_{\sing}^{\ast} \leq 3$.\\
\ \\
{\bf Theorem 16.} [O8, Theorem 26].  Let $\mathfrak{g}$ be a nonabelian, algebraic, unimodular Lie algebra such that $\trdeg_k Y(\mathfrak{g}) = i(\mathfrak{g})$.  
Suppose that $\mathfrak{g}$ admits a CP.  If $\mathfrak{g}$ is coregular then $\codim\ \mathfrak{g}_{\sing}^\ast \leq 2$.\\
\ \\
{\bf Definition 17.} We say that $\mathfrak{g}$ satisfies the Joseph-Shafrir conditions, JS for short, if $\mathfrak{g}$ is unimodular for which $p_{\mathfrak{g}}$ 
is an invariant and $\trdeg_kY(\mathfrak{g}) = i(\mathfrak{g})$.\\
Note that JS is satisfied if $\mathfrak{g}$ has no proper semi-invariants by Remark 2.\\
\ \\
The following sum rule is an extension of [OV, Proposition 1.4].\\
\ \\
{\bf Theorem 18.} [JS, Theorem 2.2]\\
Assume that $\mathfrak{g}$ satisfies JS and that $Y(\mathfrak{g})$ is freely generated by homogeneous elements $f_1,\ldots, f_r$.  Then
$$\sum\limits_{i=1}^r \deg f_i = c(\mathfrak{g}) - \deg p_{\mathfrak{g}}$$

{\bf Corollary 19.} Assume that $\mathfrak{g}$ satisfies JS and that $Y(\mathfrak{g})$ is freely generated by homogeneous elements $f_1,\ldots, f_r$.  Then
$3i(\mathfrak{g}) + 2\deg p_{\mathfrak{g}} \leq \dim\mathfrak{g} + 2 \dim Z(\mathfrak{g})$\\
Moreover, equality occurs if and only if $\deg f_i \leq 2$, $i:1,\ldots, r$.\\
\ \\
{\bf Proof.} Clearly $r = \trdeg_k Y(\mathfrak{g}) = i(\mathfrak{g})$.  We apply a similar argument as in [OV, Corollary 1.3].  The observation that $\deg f_i \geq 2$ unless $f_i \in Z(\mathfrak{g})$
combined with the preceding theorem yields:
$$\dim Z(\mathfrak{g}) + 2(i(\mathfrak{g}) - \dim Z(\mathfrak{g})) \leq \sum\limits_{i=1}^r \deg f_i = \ds\frac{1}{2} (\dim\mathfrak{g} + i(\mathfrak{g})) - \deg p_{\mathfrak{g}}$$
(and here equality occurs precisely when $\deg f_i \leq 2$ for all $i = 1,\ldots, r$)\\
$\Leftrightarrow \ -2 \dim Z(\mathfrak{g}) + 4i(\mathfrak{g}) \leq \dim \mathfrak{g} + i(\mathfrak{g}) - 2 \deg p_\mathfrak{g}$\\
$\Leftrightarrow \ 3i(\mathfrak{g}) + 2 \deg p_\mathfrak{g} \leq \dim \mathfrak{g} + 2 \dim Z(\mathfrak{g})$ \hfill $\square$ \\
\ \\
{\bf Theorem 20.} Let $\mathfrak{g}$ be unimodular, having an abelian ideal $\mathfrak{h}$ of codimension one.  Then the following are equivalent:
\begin{itemize}
\item[(1)] $\mathfrak{g}$ is coregular and $\trdeg_k Y(\mathfrak{g}) = i(\mathfrak{g})$
\item[(2)] $\mathfrak{g}$ is coregular and ad-algebraic
\item[(3)] $\dim[\mathfrak{g},\mathfrak{g}] \leq 2$
\end{itemize}

{\bf Proof.} Clearly $\mathfrak{g}$ is solvable and we may assume that $\mathfrak{g}$ is not abelian.  Then the center $Z(\mathfrak{g})$ is contained in $\mathfrak{h}$ (otherwise $\mathfrak{g} = \mathfrak{h} + Z(\mathfrak{g})$ which is abelian).  
Choose $x_0 \in \mathfrak{g} \backslash \mathfrak{h}$ and let $x_1,\ldots, x_m$, $x_{m+1},\ldots, x_n$ be a basis of $\mathfrak{h}$ such that $x_{m+1},\ldots, x_n$ is a basis  
of $Z(\mathfrak{g})$.  Then with respect to the basis $x_0, x_1, \ldots, x_m,\ldots, x_n$ of $\mathfrak{g}$ we have that $\mbox{rank}([x_i, x_j]) = 2$.  
Therefore, $i(\mathfrak{g}) = \dim \mathfrak{g} - 2$ and $c(\mathfrak{g}) = \dim\mathfrak{g}-1$.
So $\mathfrak{h}$ is a CPI of $\mathfrak{g}$.\\
We may assume that $m \geq 2$ (otherwise the result is trivial).  Put $U = \langle x_1,\ldots, x_m\rangle$.  Then
$$\mathfrak{g} = kx_0 \oplus \mathfrak{h} = kx_0 \oplus U \oplus Z(\mathfrak{g})$$
Because $\mathfrak{h}$ is abelian we see that
$$[\mathfrak{g},\mathfrak{g}] = [x_0,\mathfrak{h}] = [x_0,U]\ \ \ \mbox{and}\ \ \ C(x_0) = k x_0 \oplus Z(\mathfrak{g})$$
where $C(x_0)$ is the centralizer of $x_0$ in $\mathfrak{g}$.\\
So $\mathfrak{g} = U \oplus C(x_0)$.  Therefore $\dim C(x_0) = \dim\mathfrak{g} - \dim U$.
Next, we observe that
$$\ad~x_0 \mid_U : U \rightarrow [\mathfrak{g},\mathfrak{g}], \ \ \ x \mapsto [x_0, x]$$
is a linear bijection.  This implies that $\dim U = \dim[\mathfrak{g},\mathfrak{g}]$ and also that $[x_0, x_1] ,\ldots, [x_0, x_m]$ are linearly independent and a fortiori 
relatively prime.  Since they are the Pfaffians of the principal $2\times 2$ minors of $([x_i, x_j])$, their greatest common divisor is $p_\mathfrak{g} = 1$.  
We now proceed as follows:\\
(1) $\Rightarrow$ (3):\\
As the JS-conditions are satisfied, we may apply the preceding corollary:
$$3i(\mathfrak{g}) + 2 \deg p_\mathfrak{g} \leq \dim\mathfrak{g} + 2 \dim Z(\mathfrak{g}) < \dim\mathfrak{g} + 2 \dim C(x_0)$$
Hence, $3(\dim\mathfrak{g} - 2) < \dim\mathfrak{g} + 2(\dim\mathfrak{g} - \dim[\mathfrak{g},\mathfrak{g}])$\\
Consequently, $\dim[\mathfrak{g},\mathfrak{g}] < 3$.\\
(3) $\Rightarrow$ (2):\\
So, suppose $\dim[\mathfrak{g},\mathfrak{g}] \leq 2$.  Clearly we may assume that $\mathfrak{g}$ is indecomposable.  
This implies that $Z(\mathfrak{g}) \subset [\mathfrak{g},\mathfrak{g}]$ 
(otherwise there is a $z \in Z(\mathfrak{g})$ such that $z \notin [\mathfrak{g},\mathfrak{g}]$.  But then we could split off the abelian Lie algebra $kz$).  
Hence, $\dim Z(\mathfrak{g}) \leq 2$ and $\dim U = \dim [\mathfrak{g},\mathfrak{g}] \leq 2$. Thus $\dim \mathfrak{g}\leq 5$.  We now distinguish two cases
\begin{itemize}
\item[(i)] $\ad~x_0$ is not nilpotent\\
Decompose $\mathfrak{h}$ into the generalized weight spaces w.r.t. $\ad~x_0$:
$$\mathfrak{h} = \mathfrak{h}^o \oplus \mathfrak{h}^{\lambda_1} \oplus \ldots \oplus \mathfrak{h}^{\lambda_q}, \ \lambda_i \in k \backslash \{0\}$$
Hence, $\mathfrak{h}^{\lambda_1} \oplus \ldots \oplus \mathfrak{h}^{\lambda_q} \subset [x_0,\mathfrak{h}] = [\mathfrak{g},\mathfrak{g}]$.
Put $m_i = \dim \mathfrak{h}^{\lambda_i}$.\\
Then $m_1 + \ldots + m_q \leq \dim [\mathfrak{g},\mathfrak{g}] \leq 2$.
On the other hand, since $\mathfrak{g}$ is unimodular,
$$m_1\lambda_1 + \ldots + m_q\lambda_q = tr(\ad~x_o) = 0$$
This forces $q = 2$, $m_1 = m_2 = 1$ and $\lambda_2 = -\lambda_1$.\\
So, $\mathfrak{h}_{\lambda_1} \oplus \mathfrak{h}_{-\lambda_1} = \mathfrak{h}^{\lambda_1} \oplus \mathfrak{h}^{-\lambda_1} = [\mathfrak{g},\mathfrak{g}]$. In particular, $Z(\mathfrak{g}) = 0$\\
and $[\mathfrak{g},\mathfrak{g}] = U = \mathfrak{h}$. Choose nonzero $y_1 \in \mathfrak{h}_{\lambda_1}$, $y_2 \in \mathfrak{h}_{-\lambda_1}$ and put $y_0 = (1/\lambda_1)x_0$.\\
Then, $y_0$, $y_1$, $y_2$ is a basis for $\mathfrak{g}$ with nonzero brackets: $[y_0, y_1] = y_1$, $[y_0, y_2] = -y_2$. Clearly $\mathfrak{g}$ is algebraic and also coregular since $Y(\mathfrak{g}) = k[y_1y_2]$.
\item[(ii)] $\ad~x_0$ is nilpotent.\\
In this case $\mathfrak{g}$ is nilpotent and thus algebraic.  Consulting [D2, O7], the following are the only indecomposable nilpotent Lie algebras of dimension at most 5 having a commutative
ideal of codimension one and such that $\dim [\mathfrak{g},\mathfrak{g}] \leq 2 : \mathfrak{g}_3, \mathfrak{g}_4, \mathfrak{g}_{5,2}$, which are all coregular.
\end{itemize}
(2) $\Rightarrow$ (1): Denote by $Q(Y(\mathfrak{g}))$ the quotient field of $Y(\mathfrak{g})$.\\
It suffices to show that $Q(Y(\mathfrak{g})) = R(\mathfrak{g})^{\mathfrak{g}}(\bullet)$, because then we obtain at once that
$$\trdeg_kY(\mathfrak{g}) = \trdeg_kR(\mathfrak{g})^\mathfrak{g}= i(\mathfrak{g})$$
by Theorem 1 since $\mathfrak{g}$ is ad-algebraic.\\
Let $\ad~x_o = S + N$ be the Jordan decomposition of $\ad~x_0$, with $S$ and $N$ its semi-simple and nilpotent components. As $\mathfrak{g}$ is ad-algebraic we can find $s, y \in \mathfrak{g}$ such that $S = \ad~s$ and $N = \ad~y$.\\
We distinguish 2 cases:
\begin{itemize}
\item[a)] $\ad~s(\mathfrak{h}) \neq 0$, i.e. $s \notin \mathfrak{h}$.  Then we replace $x_0$ by a suitable nonzero scalar multiple of $s$, which is diagonalizable with integer eigenvalues having zero sum.  By the same argument
as in the proof of [O8, Example 28] we obtain $(\bullet)$.
\item[b)] $\ad~s(\mathfrak{h}) = 0$.  Then $\ad~y(\mathfrak{h}) \neq 0$ (since $\ad~x_0 (\mathfrak{h}) \neq 0)$, i.e. $y \notin \mathfrak{h}$.  
Then we replace $x_0$ by $y$.  It follows that $\mathfrak{g}$ is nilpotent for which $(\bullet)$ is
well known (since $\mathfrak{g}$ has no proper semi-invariants).\hfill $\square$
\end{itemize}
{\bf Remark 21.} A more direct approach for the implication (2) $\Rightarrow$ (3) goes as follows.\\
By assumption $\mathfrak{g}$ is nonabelian, unimodular, coregular and ad-algebraic.  As above it then also satisfies $\trdeg_kY(\mathfrak{g}) = i(\mathfrak{g})$.  Moreover $\mathfrak{h}$ is a CP of $\mathfrak{g}$.
Next, we observe that 
$$\mathfrak{g}_{\sing}^\ast = \{\xi \in \mathfrak{g}^\ast \mid \xi ([x_0, x_i]) = 0,\ i=1,\ldots, m\}$$
Then, $\dim [\mathfrak{g},\mathfrak{g}] = m = \codim~\mathfrak{g}_{\sing}^{\ast} \leq 2$ by Theorem 16.\hfill $\square$\\
\ \\
{\bf Remark 22.} In the list of all indecomposable nilpotent Lie algebras of dimension at most seven [O7,O8] there are only 6 Lie algebras with an 
abelian ideal of codimension one and with $\dim [\mathfrak{g},\mathfrak{g}] > 2$, namely 8, 25, 156, 157, 158, 159.  None of these
is coregular as predicted by Theorem 20.\\
Examples 24-26 show that none of the conditions such as unimodular, algebraic and $\trdeg_kY(\mathfrak{g}) = i(\mathfrak{g})$ can be removed from Theorem 18 and Corollary
19.  In all, except for Example 26, $\mathfrak{h} = \langle x_2, x_3, x_4\rangle$ is an abelian ideal of codimension one.\\
\ \\
{\bf Example 23.} Let $\mathfrak{g}$ be the solvable Lie algebra with basis $x_1$, $x_2$, $x_3$, $x_4$ and with nonzero brackets
$$[x_1,x_2] = x_2, \ [x_1, x_3] = x_3, \ [x_1,x_4] = -2x_4$$
Clearly, $\mathfrak{g}$ is unimodular, algebraic and $\dim [\mathfrak{g},\mathfrak{g}] = 3 > 2$.  Hence $\mathfrak{g}$ is not coregular by Theorem 20 (see also [JS, 8.4] and
[O8, Example 28]).  This can also be seen directly.\\
Indeed, $Y(\mathfrak{g}) = k[f_1, f_2, f_3]$ where $f_1 = x_2^2x_4$, $f_2 = x_3^2 x_4$, $f_3 = x_2x_3x_4$ with $f_1f_2 = f_3^2$.\\
In particular, $Y(\mathfrak{g})$ is not factorial.
Note that $Sy(\mathfrak{g}) = k[x_2, x_3, x_4]$, which is polynomial, and $R(\mathfrak{g})^\mathfrak{g}= k(f_1, f_3)$.  So $\trdeg_k(Y(\mathfrak{g}) = 2 = i(\mathfrak{g})$ while
$p_\mathfrak{g} =1$.  Thus JS is satisfied.  Finally, $\codim \ \mathfrak{g}_{\sing}^{\ast} = 3$.\\
\ \\
{\bf Example 24.}\\
Let $\mathfrak{g}$ be the Lie algebra with basis $x_1$, $x_2$, $x_3$, $x_4$ and nonzero brackets
$$[x_1,x_2] = x_2, \ [x_1, x_3] = x_3, \ [x_1, x_4] = -x_4$$
Clearly, $\mathfrak{g}$ is algebraic, but not unimodular.  Also, $p_\mathfrak{g} = 1$ and $Z(\mathfrak{g}) = 0$.  Put $f_1 = x_2x_4$ and $f_2 = x_3x_4$. Then $Y(\mathfrak{g}) = k[f_1,f_2]$,
so $\mathfrak{g}$ is coregular and $\trdeg_kY(\mathfrak{g}) = 2 = i(\mathfrak{g})$.  Moreover,
$$Sy(\mathfrak{g}) = k[x_2, x_3, x_4] \ \mbox{and}\ R(\mathfrak{g})^{\mathfrak{g}} = k(f_1,f_2)$$
However,
$$\deg f_1 + \deg f_2 = 4 > 3 = c(\mathfrak{g}) - \deg p_\mathfrak{g}$$
$$3i(\mathfrak{g}) + 2\deg p_\mathfrak{g} = 6 > 4 = \dim \mathfrak{g} + 2 \dim Z(\mathfrak{g})$$
and $\dim [\mathfrak{g},\mathfrak{g}] = 3 >2$.\\
\ \\
{\bf Example 25.}\\
Let $\mathfrak{g}$ be the Lie algebra with basis $x_1$, $x_2$, $x_3$, $x_4$ and nonzero brackets
$$[x_1, x_2] = x_2 + x_3,\ [x_1, x_3] = x_3, \ [x_1, x_4] = -2x_4$$
$\mathfrak{g}$ is unimodular, but not algebraic (not even almost algebraic).\\
Again, $p_\mathfrak{g} = 1$ and $Z(\mathfrak{g}) = 0$. $\mathfrak{g}$ is coregular since $Y(\mathfrak{g}) = k[x_3^2 x_4]$.\\
However, $\trdeg_kY(\mathfrak{g}) = 1 < 2 = i(\mathfrak{g})$.  Clearly,
$$Sy(\mathfrak{g}) = k[x_3, x_4]\ \mbox{and}\ R(\mathfrak{g})^\mathfrak{g} = k(x_3^2x_4)$$
In particular, $j(\mathfrak{g}) = 1 < i(\mathfrak{g})$.  Also,
$$3i(\mathfrak{g}) + 2\deg p_\mathfrak{g} = 6 > 4 = \dim \mathfrak{g} + 2 \dim Z(\mathfrak{g})$$
and $\dim [\mathfrak{g},\mathfrak{g}] = 3 > 2$.  Finally, we notice that $\mathfrak{g}_\Lambda = \langle x_2, x_3, x_4\rangle = F(\mathfrak{g})$ and $Sy(\mathfrak{g}) \neq k[x_2,x_3,x_4] = Y(\mathfrak{g}_\Lambda)$.\\
\ \\
{\bf Example 26.}\\
Consider the 9-dimensional solvable Lie algebra $\mathfrak{g}$ with basis $x_0, x_1,\ldots, x_8$ and nonzero brackets\\
$[x_0, x_1] = 5x_1$, $[x_0, x_2] = 10x_2$, $[x_0, x_3] = -13x_3$, $[x_0, x_4] = -8x_4$,\\
$[x_0, x_5] = -3x_5$, $[x_0, x_6] = 2x_6$, $[x_0, x_7] = 7x_7$, $[x_1, x_3] = x_4$,\\
$[x_1, x_4] = x_5$, $[x_1, x_5] = x_6$, $[x_1,x_6] = x_7$, $[x_2, x_3] = x_5$, $[x_2, x_4] = x_6$, $[x_2, x_5] = x_7$.\\
Then, $\mathfrak{g}$ is algebraic and unimodular with $\codim \ \mathfrak{g}_{\sing}^{\ast} = 3$.  In particular, $p_\mathfrak{g} = 1$.\\
$\mathfrak{g}$ is coregular since $Y(\mathfrak{g}) = k[x_8]$, but $\trdeg_k Y(\mathfrak{g}) = 1 < 3 = i(\mathfrak{g})$.\\
Note that $c(\mathfrak{g}) = 6$ and $\deg x_8 = 1 < 6 = c(\mathfrak{g}) - \deg p_\mathfrak{g}$.\\
So the sum rule fails in these circumstances.  Furthermore, $\mathfrak{g}_\Lambda = \langle x_1,\ldots, x_8\rangle$.  By (5) of Theorem 4 and [DDV, p.323]
$$Sy(\mathfrak{g}) = Y(\mathfrak{g}_\Lambda) = k[x_7, x_8, f, g, h]$$
where
\begin{eqnarray*}
f &=& 3x_4 x_7^2 - 3x_5 x_6 x_7 + x_6^3\\
g &=& 4x_3 x_7^2 - 2x_5^2x_7^2 - 4x_4 x_6 x_7^2 + 4x_5x_6^2 x_7 - x_6^4\\
h &=& (f^4 + g^3)/x_7^3
\end{eqnarray*}
Hence $Sy(\mathfrak{g}) $ is not polynomial.  Note that $\mathfrak{g}_\Lambda$ is isomorphic to a central extension of the nilpotent Lie algebra with number 152 of [O8, p.109].  
Finally, $R(\mathfrak{g})^\mathfrak{g} = k(x_8, f^4 g^{-3}, fgx_7^{-2})$ and $F(\mathfrak{g}) = \langle x_3, x_4, x_5, x_6, x_7, x_8\rangle$ which is a CPI of $\mathfrak{g}$.\\
\ \\
{\bf Remark 27.} The preceding example is a central extension of example (58) of [DDV, p.322], which turned out to be a counterexample to Bolsinov's completeness criterion
for Mishchenko-Fomenko subalgebras [O8, Counterexample 20]. Inspired by this, Bolsinov obtained an interesting and useful adaptation of his original criterion by considering an alternative
definition for Mishchenko-Fomenko subalgebras [Bo].  See also [JS, Theorem 7.2].\\
\ \\
Let $\mathfrak{g}$ be a semi-simple Lie algebra, $B$ a Borel subalgebra of $\mathfrak{g}$.  
Then it is well known that the nilradical of $B$ is coregular [J2, 4.7], see also Corollary 32.  An example by A. Hersant shows that a similar result does not hold in general if we 
replace $B$ by an arbitrary parabolic subalgebra of $\mathfrak{g} $ [J2, 8.5]. \\
We will now give a short proof of an extension of this example.\\
\ \\
{\bf Proposition 28.}\\
Let $N$ be the nilradical of the parabolic subalgebra $P$ of type $(1,1,n-2)$ inside $sl(n)$, with $n \geq 3$.  Then
$$N \ \mbox{is coregular}\ \ \ \Leftrightarrow \ \ \ n \leq 4 \ \ \ \Leftrightarrow \ \ \ i(N) \leq 3$$

{\bf Proof.}\\
Let $(E_{ij})$, $i,j=1,\ldots, n$, be the standard basis for $gl(n)$.  Then
$$\{E_{12}, E_{13},\ldots, E_{1n}; E_{23},\ldots, E_{2n}\}$$
is a basis for $N$ (so $\dim N = 2n-3$), with nonzero brackets
$$[E_{12}, E_{23}] = E_{13},\ [E_{12}, E_{24}] = E_{14},\ldots, [E_{12}, E_{2n}] = E_{1n}$$
Clearly, $[N,N] = \langle E_{13}, E_{14},\ldots,E_{1n}\rangle$ and $\dim [N,N] = n-2.$\\
$N$ admits an abelian ideal of codimension one, namely $H = \langle E_{13},\ldots, E_{1n}; E_{23},\ldots, E_{2n}\rangle$.\\
$N$, being nilpotent, is ad-algebraic and $i(N) = 2n-5$.\\
By Theorem 20:
$$N\ \mbox{is coregular}\ \ \ \Leftrightarrow \ \ \ \dim [N,N] \leq 2 \ \ \ \Leftrightarrow \ \ \ n \leq 4 \ \ \ \Leftrightarrow \ \ \ i(N) \leq 3$$
\hfill $\square$\\
\ \\
{\bf 3.2 Sufficient conditions for polynomiality}\\
We exhibit some methods which will be used in sections 4 and 5.  The first one is very efficient for proving coregularity, provided one has candidates for the generating
invariants.  It is an extension of [PPY, Theorem 1.2]. See also [Pa2, Theorem 1.2].\\
\ \\
{\bf Theorem 29.} [JS, 5.7], [Sh]\\
Assume that $f_1, \ldots, f_r \in Y(\mathfrak{g})$, $r = i(\mathfrak{g})$, are algebraically independent homogeneous invariants such that
$$\sum\limits_{i=1}^r \deg f_i \leq c(\mathfrak{g})) - \deg p_\mathfrak{g}$$
Then, equality holds and $Y(\mathfrak{g}) = k[f_1,\ldots, f_r]$.  In particular, $\trdeg_kY(\mathfrak{g}) = i(\mathfrak{g})$.\\
\ \\
{\bf Theorem 30.} [O7, Theorem 3.5] (The Frobenius method).\\
Let $\mathfrak{g}$ be a finite dimensional Lie algebra over $k$.  Assume that there exists a torus $T \subset \Der \mathfrak{g}$ (i.e. an abelian subalgebra consisting of semi-simple derivations 
of $\mathfrak{g}$) such that the semi-direct product $L = T \oplus \mathfrak{g}$ is Frobenius.  Let $f_1,\ldots, f_r$ be the irreducible factors of $p_L$ (equivalently of 
$\Delta(L)$).  Then the following hold:
\begin{itemize}
\item[(1)] $Sy(\mathfrak{g}) = Sy(L) = Y(L_\Lambda) = k [f_1, \ldots, f_r]$, a polynomial algebra.
\item[(2)] $\dim T = i(\mathfrak{g})$ and $r = i(L_\Lambda) = \dim L - \dim L_\Lambda$.
\item[(3)] $\Lambda(\mathfrak{g}) = \{\lambda\mid_\mathfrak{g} \mid \lambda \in \Lambda (L)\}$ and $\mathfrak{g}_\Lambda = \mathfrak{g} \cap L_\Lambda$
\item[(4)] If $\mathfrak{g}$ has no proper semi-invariants (i.e. $\mathfrak{g} = \mathfrak{g}_\Lambda$) then $\mathfrak{g} = L_\Lambda$ and
$$Y(\mathfrak{g}) = k[f_1,\ldots, f_r]$$
\end{itemize}

{\bf Remark 31.}  Although this method does not always work, it has some significant advantages.  First of all it is relatively simple: it comes down to showing that the 
determinant $\Delta (L)$
of the structure matrix of $L$ is not zero.  In addition, there is no need to have prior knowledge of candidates for the generating (semi-) invariants.  In fact, we get them as a bonus since they are precisely
the irreducible factors of the determinant above (or equivalently of $p_L$).  This method works rather well if $\mathfrak{g}$ is nilpotent.  It will also be useful in sections 4 and 5.  In [O7,O8] the Poisson
center has been determined explicitly for the 159 cases of the indecomposable nilpotent Lie algebras of dimension at most seven (here a family is counted as one Lie algebra).  It turns out that 132 of them are coregular.  Among the latter, 67 Lie algebras were
treated successfully with this method [O7, 5].\\
\ \\
{\bf Corollary 32.} See also [J2, 4.7].  Let $\mathfrak{g}$ be a simple Lie algebra with triangular decomposition $\mathfrak{g} = N^- \oplus H \oplus N$. Then the nilradical $N$ of the Borel subalgebra
$B = H \oplus N$ is coregular.\\
\ \\
{\bf Proof.}  There exists a torus $T \subset \ad_N H \subset \Der N$ such that the semi-direct product $L = T\oplus N$ is Frobenius.
Hence $Y(N)$ is polynomial by (4) of Theorem 30.  Indeed, in case $\mathfrak{g}$ is not of type
$A_n$, $n \geq 2$; $D_{2t+1}$, $t\geq 2$; or $E_6$ then it suffices to take $T = \ad_NH$ because then $T \oplus N = B$, which is Frobenius (for more details see [O7, Corollary 3.6]).
The existence of $T$ if $\mathfrak{g}$ is of type $A_n$ is easy to verify.  The remaining cases were done by Rupert Yu (unpublished). \hfill $\square$\\
\ \\
{\bf Question 33.} (Rupert Yu)\\
Suppose $\mathfrak{g}$ is a Lie algebra for which there exists a derivation $d \in \Der \mathfrak{g}$ such that the semidirect product $L = kd \oplus \mathfrak{g}$ is Frobenius.
Does this imply that $Sy(\mathfrak{g}) = Sy(L)$ ?\\
\ \\
We know this is true if $d$ is diagonalizable by Theorem 30.  However the following is a counterexample for the general case.\\
\ \\
{\bf Example 34.}  Let $\mathfrak{g}$ be the 5-dimensional Lie algebra with basis $x_1,\ldots, x_5$ and nonzero brackets: $[x_1, x_3] = x_3 - x_4,\ [x_1, x_4] = x_4,  
[x_1, x_5] = x_5, \ [x_2, x_3] = x_5$.\\
$\mathfrak{g}$ is solvable of index one, but it is not almost algebraic.  One verifies that 
$$Y(\mathfrak{g}) = k, \ Sy(\mathfrak{g}) = k[x_4, x_5] \ \mbox{and}\ R(\mathfrak{g})^{\mathfrak{g}} = k (x_4/x_5)$$
Note that $\trdeg_kY(\mathfrak{g}) = 0 < 1 = i(\mathfrak{g})$.  Also, $j(\mathfrak{g}) = i(\mathfrak{g})$.\\
Next we take the derivation $d \in \Der~\mathfrak{g}$ given by
$$d(x_1) = - x_2, \ d(x_2) = 0, \ d(x_3) = x_4, \ d(x_4) = x_5, \ d(x_5) = 0$$
Clearly $d$ is nilpotent.  Consider $L = kd\oplus \mathfrak{g}$.  Then, $\Delta(L) = x_5^6 \neq 0$.\\
Hence $L$ is Frobenius and $Sy(L) = k[x_5]$ (by Theorem 30),  which does not coincide with $Sy(\mathfrak{g})$.\\
\ \\
{\bf Definition 35.} $Y(\mathfrak{g})$ is said to be saturated if for some nonzero $u$, $v \in S(\mathfrak{g})$, $uv \in Y(\mathfrak{g})$ implies that so are $u$ and $v$.  
In particular, $Y(\mathfrak{g})$ is factorial.  Note that the condition $u,v \in S(\mathfrak{g})$ may be replaced by $u, v \in Sy(\mathfrak{g})$ 
because $uv \in Y(\mathfrak{g})$ implies that $u$ and $v$
are semi-invariants and thus belong to $Sy(\mathfrak{g})$.\\
We now recall when $Sy(\mathfrak{g})$ is a polynomial ring over $Y(\mathfrak{g})$. Clearly, in order for this to happen $Y(\mathfrak{g})$ must be saturated.\\
\ \\
{\bf Theorem 36.} [DNOW, Theorem 6]\\
Assume that
\begin{itemize}
\item[(i)] $R(\mathfrak{g})^\mathfrak{g} = Q(Y(\mathfrak{g}))$, the quotient field of $Y(\mathfrak{g})$.
\item[(ii)] $Y(\mathfrak{g})$ is saturated.
\end{itemize}
Then the following hold:
\begin{itemize}
\item[(1)] $S(\mathfrak{g})$ has at most a finite number of irreducible proper semi-invariants $v_1,\ldots, v_t$.\\
Let $\lambda_1,\ldots, \lambda_t \in \Lambda(\mathfrak{g})$ be their weights.
\item[(2)] $S(\mathfrak{g})_{\lambda_i} = Y(\mathfrak{g}) v_i$, $i = 1,\ldots, t$.
\item[(3)] Each $v_i$ is a semi-invariant for all derivations $d \in \Der \ \mathfrak{g}$.
\item[(4)] $Sy(\mathfrak{g}) = Y(\mathfrak{g})[v_1, \ldots, v_t]$, a polynomial ring over $Y(\mathfrak{g})$.  In particular, if $Y(\mathfrak{g})$ is polynomial over $k$, then
the same holds for $Sy(\mathfrak{g})$. [We don't know if the converse holds. It is a special case of the Zariski cancellation problem.  For the general question Susumu Oda claims to have a proof [Oda], but some experts are skeptical]
\end{itemize}
We now look at a special case of Theorem 36.\\
\ \\
{\bf Theorem 37.} [DNOW, Proposition 16 and Theorem 18]  Let $x_1,\ldots, x_s$, $x_{s+1},\ldots, x_n$ be a basis such that $x_1,\ldots, x_s$ is a basis of $Z(\mathfrak{g})$.
Let $E_1,\ldots, E_m$ be a basis of the algebraic hull $H$ of $\ad~\mathfrak{g}$.  Then the following conditions are equivalent:
\begin{itemize}
\item[(1)] $j(\mathfrak{g}) = \dim Z(\mathfrak{g})$ (i.e. $\mathfrak{g}[\xi] = Z(\mathfrak{g}))$ for some $\xi \in \mathfrak{g}^{\ast}$)
\item[(2)] $R(\mathfrak{g})^\mathfrak{g} = R(Z(\mathfrak{g}))$
\item[(3)] $Z(D(\mathfrak{g})) = D(Z(\mathfrak{g})) = k (x_1,\ldots, x_s)$, a rational extension of $k$.
\item[(4)] The localization $U(\mathfrak{g})_S$, where $S = U(Z(\mathfrak{g})) \backslash \{0\}$, is primitive.
\end{itemize}
Moreover, these conditions imply that:
\begin{itemize}
\item[(a)] $Y(\mathfrak{g}) = S(Z(\mathfrak{g})) = k[x_1,\ldots, x_s]$, which is saturated.
\item[(b)] $S(\mathfrak{g})$ admits at most a finite number of irreducible, proper semi-invariants $v_1, \ldots, v_t$.
\item[(c)] $Sy(\mathfrak{g}) = Y(\mathfrak{g}) [v_1,\ldots, v_t] = k[x_1,\ldots, x_s, v_1,\ldots, v_t]$, a polynomial algebra over $k$.
\item[(d)] $v_1,\ldots, v_t$ are precisely the irreducible factors, not in $Y(\mathfrak{g})$, of $p'_{\mathfrak{g}} \in S(\mathfrak{g})$, the latter being the greatest common
divisor of the $r \times r$ minors of the $m \times n$ matrix $(E_ix_j)$, where
$$r = \rank (E_ix_j) = \dim \mathfrak{g} - \dim Z(\mathfrak{g})$$
\end{itemize}

{\bf Remark 38.} Suppose $\mathfrak{g}$ is square integrable, i.e. $i(\mathfrak{g}) = \dim Z(\mathfrak{g})$ which forces $j(\mathfrak{g}) = \dim Z(\mathfrak{g})$.  So, the above
conditions are satisfied and we may replace the matrix $(E_ix_j)$ by the structure matrix $([x_i, x_j])$ of $\mathfrak{g}$.  Consequently $v_i, \ldots, v_t$ are then precisely the
irreducible factors, not in $Y(\mathfrak{g})$, of $p_{\mathfrak{g}}$.\\
\ \\
{\bf Example 39.} Let $\mathfrak{g}$ be the 4-dimensional Lie algebra with basis $x_1$, $x_2$, $x_3$, $x_4$ and nonzero brackets $[x_1,x_2] = x_2 + x_3$, $[x_1, x_3] = x_4$. \\
Clearly,
$i(\mathfrak{g}) = 2$, $Z(\mathfrak{g}) = \langle x_4\rangle$ and $p_\mathfrak{g} = 1$.\\
$\dim Z(\mathfrak{g}) = 1 < i(\mathfrak{g})$, so $\mathfrak{g}$ is not square integrable.  Obviously, in this situation $p_\mathfrak{g}$ is useless in order to compute the remaining semi-invariants.
Next, we introduce $E_1, E_2 \in \Der \ \mathfrak{g}$ as follows:\\
$E_1 (x_1) = 0$, $E_1(x_2) = x_2 + x_2 + x_4$, $E_1(x_3) = E_1 (x_4) = 0$\\
$E_2(x_1) = 0$, $E_2(x_2) = -x_4$, $E_2(x_3) = x_4$, $E_2(x_4) = 0$.\\
In fact, $E_1$ and $E_2$ are the semi-simple and nilpotent components of $\ad~x_1$.\\
Hence they belong to the algebraic hull $H$ of $\ad~\mathfrak{g}$.  One verifies that 
$$E_1, E_2, E_3 = \ad~x_2, E_4 = ad \ x_3$$
form a basis of $H$.  We now observe the matrix $(E_i x_j)$:

\begin{center}
\begin{tabular}{l|cccc}
 &$x_1$ &$x_2$ &$x_3$ &$x_4$\\
 \hline
 $E_1$ &0 &$x_2 + x_3 + x_4$ &0 &0\\
 $E_2$ &0 &$-x_4$ &$x_4$ &0\\
 $E_3$ &$-x_2-x_3$ &0 &0 &0\\
 $E_4$ &$-x_4$ &0 &0 &0
 \end{tabular}
 \end{center}
 which is of rank 3.  By Theorem 1
 $$\trdeg_kR(\mathfrak{g})^{\mathfrak{g}} = j(\mathfrak{g}) = \dim \mathfrak{g} - \rank~(E_ix_j) = 1.$$
 Since $j(\mathfrak{g}) = 1 = \dim Z(\mathfrak{g})$ we can apply Theorem 37.  Hence,
 $$Y(\mathfrak{g}) = k[x_4] \ \mbox{and}\ R(\mathfrak{g})^{\mathfrak{g}} = k(x_4)$$
 Furthermore, the matrix $(E_ix_j)$ has only 2 nonzero $3\times 3$ minors, namely
$$ -x_4(x_2 + x_3) (x_2 + x_3 + x_4)\ \mbox{and}\ -x_4^2(x_2 + x_3 + x_4)$$
Their greatest common divisor is $p'_{\mathfrak{g}} = x_4(x_2 + x_3 + x_4)$.
By (c) of Theorem 37 we may conclude that
$$Sy(\mathfrak{g}) = Y(\mathfrak{g}) [x_2 + x_3 + x_4] = k[x_4, x_2 + x_3 + x_4]$$

{\bf 4. Coregularity for Lie algebras with index at most two}\\
{\bf Motivation:} Due to Corollary 19 there are not many nonabelian, coregular Lie algebras with a large index.  
This is especially true if the index is maximal (Proposition 40 and Corollary 41).  On the other hand,
we will encounter quite a few Lie algebras for which the coregularity implies that their index is at most two 
(Proposition 42, Proposition 50, Theorem 51, subsection 5.1). \\
\ \\
{\bf Proposition 40}\\
Assume that $\mathfrak{g}$ is an indecomposable Lie algebra which satisfies JS and for which $i(\mathfrak{g}) = \dim \mathfrak{g}-2$.  If $\mathfrak{g}$ is coregular then either $\mathfrak{g} = sl(2,k)$ or $\mathfrak{g}$
is solvable with $\dim \mathfrak{g} \leq 6$.\\
\ \\
{\bf Proof.} First we notice that $Z(\mathfrak{g}) \subset [\mathfrak{g},\mathfrak{g}]$ as $\mathfrak{g}$ is indecomposable.  Next we claim that 
$\dim Z(\mathfrak{g}) \leq \ds\frac{1}{2} \dim \mathfrak{g}$.\\
Indeed, take $\xi \in \mathfrak{g}_{\reg}^\ast$.  Then the stabilizer $\mathfrak{g}(\xi)$ is abelian [D6, 1.11.7] of dimension $i(\mathfrak{g}) = \dim \mathfrak{g}-2$ and $Z(\mathfrak{g}) \subset \mathfrak{g}(\xi)$.  There exists a basis\\
$x_1, x_2, x_3, \ldots, x_p, x_{p+1},\ldots, x_n$ of $\mathfrak{g}$ such that $x_3,\ldots, x_n$ is a basis of $\mathfrak{g}(\xi)$ and $x_{p+1},\ldots, x_n$ is a basis of $Z(\mathfrak{g})$.  It suffices to show that $\dim [\mathfrak{g},\mathfrak{g}] \leq p$, because then
$$2 \dim Z(\mathfrak{g}) \leq \dim [\mathfrak{g},\mathfrak{g}] + \dim Z(\mathfrak{g}) \leq p + \dim Z(\mathfrak{g}) = \dim \mathfrak{g}$$
Clearly, $[x_1,x_2] \neq 0$ and the structure matrix $M = ([x_i,x_j])_{1 \leq i,j \leq n}$ of $\mathfrak{g}$ has rank $r = n-i(\mathfrak{g}) = 2$.\\
We may assume that $[x_1, x_3] \neq 0$ and that $p > 3$ (if $p = 3$ then $\dim [\mathfrak{g},\mathfrak{g}] = \dim \langle [x_1, x_2], [x_1, x_3], [x_2,x_3]\rangle \leq 3 = p$).  This implies that the following submatrix $A$ of $M$ has rank one (otherwise rank $M = 4$)
$$A = \left(\begin{array}{l}
[x_1,x_3] \ldots [x_1,x_p]\\[1ex]
[x_2,x_3] \ldots [x_2,x_p]\end{array}\right)$$
Using the fact that its nonzero entries have degree one, it is not difficult to see that we have to consider the following two cases:
\begin{itemize}
\item[(1)] Each column of $A$ is a scalar multiple of the first one.  Then,
$$\dim [\mathfrak{g},\mathfrak{g}] = \dim \langle [x_1,x_2], [x_1,x_3],[x_2, x_3]\rangle \leq 3 < p$$
\item[(2)] The second row of $A$ is a scalar multiple of the first one.  Then,
$$\dim [\mathfrak{g},\mathfrak{g}] = \dim \langle [x_1,x_2], [x_1,x_3], \ldots, [x_1, x_p]\rangle \leq p$$
This establishes the claim.
\end{itemize}
Application of Corollary 19 gives us:
$$3(\dim \mathfrak{g} -2) = 3i(\mathfrak{g}) \leq \dim \mathfrak{g} + 2\dim Z(\mathfrak{g}) \leq 2 \dim \mathfrak{g}$$
Consequently, $\dim \mathfrak{g} \leq 6$.  Hence $\mathfrak{g}$ is solvable or $\mathfrak{g} = sl(2,k)$ (otherwise $i(\mathfrak{g}) < \dim \mathfrak{g} - 2$ by [AOV2, pp. 554-559] 
or by subsection 5.1). \hfill $\square$\\
\ \\
{\bf Corollary 41.} Suppose $\mathfrak{g}$ is an indecomposable nilpotent Lie algebra with $i(\mathfrak{g}) = \dim \mathfrak{g} - 2$.  If $\mathfrak{g}$ is coregular then by the 
above and [O7, 5] $\mathfrak{g}$ is isomorphic to one of the following:
$$\mathfrak{g}_3, \ \ \mathfrak{g}_4, \ \ \mathfrak{g}_{5,2}, \ \ \mathfrak{g}_{5,4}, \ \ \mathfrak{g}_{6,3}$$

{\bf Proposition 42.} Assume that $\mathfrak{g}$ satisfies JS with $\dim Z(\mathfrak{g}) \leq 1$. If $\mathfrak{g}$ is coregular then $i(\mathfrak{g})$ is 1 or 2 in each of the following
cases:
\begin{itemize}
\item[(1)] $7 \neq \dim \mathfrak{g} \leq 8$
\item[(2)] $\dim \mathfrak{g} = 7$ and $\mathfrak{g}$ is singular
\item[(3)] $\dim \mathfrak{g} = 9$ or 10 and $\deg p_\mathfrak{g} \geq 2$
\end{itemize}

{\bf Proof.} Again the main tool will be Corollary 19.  Being unimodular, $\mathfrak{g}$ is not Frobenius [O3, Theorem 3.3].  Hence, $i(\mathfrak{g}) \geq 1$.
\begin{itemize}
\item[(1)] First we suppose $\dim \mathfrak{g}$ (and hence also $i(\mathfrak{g})$) is even.  Then
$$3i(\mathfrak{g}) \leq 3i(\mathfrak{g}) + 2\deg p_\mathfrak{g} \leq \dim \mathfrak{g} + 2 \dim Z(\mathfrak{g}) \leq 10$$
implies that $i(\mathfrak{g}) = 2$.\\
On the other hand, if $\dim \mathfrak{g}$ (and hence also $i(\mathfrak{g})$) is odd, then
$$3i(\mathfrak{g}) \leq 3i(\mathfrak{g}) + 2\deg p_\mathfrak{g} \leq \dim \mathfrak{g} + 2 \dim Z(\mathfrak{g}) \leq 7$$
forces $i(\mathfrak{g}) = 1$.
\item[(2)] By assumption $\dim \mathfrak{g} = 7$ and $\deg p_\mathfrak{g} \geq 1$ since $\mathfrak{g}$ is singular.  Hence,
$$3i(\mathfrak{g}) < 3i(\mathfrak{g}) + 2\deg p_\mathfrak{g} \leq \dim \mathfrak{g} + 2 \dim Z(\mathfrak{g}) \leq 9$$
and thus $i(\mathfrak{g}) = 1$.
\item[(3)]follows at once from $\deg p_\mathfrak{g} \geq 2$ and
$$3i(\mathfrak{g}) + 2\deg p_\mathfrak{g} \leq \dim \mathfrak{g} + 2 \dim Z(\mathfrak{g}) \leq 12$$
\end{itemize}
\hfill $\square$\\
\ \\
{\bf Proposition 43.} Let $\mathfrak{g}$ be a Lie algebra with $j(\mathfrak{g}) \leq 1$ (as it is when $i(\mathfrak{g}) \leq 1$).
\begin{itemize}
\item[(a)] If $j(\mathfrak{g}) = 0$ then $Z(\mathfrak{g}) = 0$. By Theorem 37
$$Y(\mathfrak{g}) = k,\ R(\mathfrak{g})^\mathfrak{g} = k \ \ \ \mbox{and}\ \ \ Sy(\mathfrak{g})\ \mbox{is polynomial.}$$
\item[(b)] Let $j(\mathfrak{g}) = 1$ (in particular $\trdeg_kR(\mathfrak{g})^\mathfrak{g} = 1$) and suppose $Y(\mathfrak{g}) \neq k$.  Choose a homogeneous
element $v \in Y(\mathfrak{g}) \backslash k$ of smallest degree.  Then by [O3, Lemma 3.8], $Y(\mathfrak{g}) = k[v]$ and $R(\mathfrak{g})^\mathfrak{g} = k(v)$.  Now, assume in addition that
$Y(\mathfrak{g})$ is saturated.  Then $v$ is irreducible and by Theorem 36 there are irreducible semi-invariants $v_1,\ldots, v_t$ in $S(\mathfrak{g})$ such that 
$$Sy(\mathfrak{g}) = Y(\mathfrak{g}) [v_1,\ldots, v_t] = k[v, v_1,\ldots, v_t]$$
which is a polynomial algebra. \hfill $\square$\\
\end{itemize}

Now we need to recall a special case of a result by Dixmier [D1, p. 333]:\\
\ \\
{\bf Theorem 44.} Let $\mathfrak{g}$ be a nilpotent Lie algebra and let
$$0 = \mathfrak{g}_0 \subset \mathfrak{g}_1 \ldots \subset \mathfrak{g}_n = \mathfrak{g}$$
be a sequence of ideals of $\mathfrak{g}$ such that for each $j : 1,\ldots,n$, $\dim \mathfrak{g}_j = j$ and $[\mathfrak{g},\mathfrak{g}_j] \subset \mathfrak{g}_{j-1}$.  Choose $x_j \in \mathfrak{g}_j \backslash \mathfrak{g}_{j-1}$.  Suppose $j_1 < j_2 < \ldots < j_r$ are the indices $j \geq 1$ such that
$$S(\mathfrak{g}_{j-1}) \cap Y(\mathfrak{g}) \stackrel{\subset}{\neq} S(\mathfrak{g}_j) \cap Y(\mathfrak{g})$$
\begin{itemize}
\item[(1)] Then for each such $j$ there is a nonzero element $b_j \in S(\mathfrak{g}_{j-1}) \cap Y(\mathfrak{g})$ and $c_j \in S(\mathfrak{g}_{j-1})$ such that $a_j = b_j x_j + c_j \in S(\mathfrak{g}_j) \cap Y(\mathfrak{g})$.
\end{itemize}
In (2), (3), (4) $a_j$, $b_j$, $c_j$ are chosen to satisfy  (1).
\begin{itemize}
\item[(2)] $Y(\mathfrak{g}) \subset k[a_{j_1},\ldots, a_{j_r}, b_{j_1}^{-1},\ldots, b_{j_r}^{-1}]$
\item[(3)] $R(\mathfrak{g})^{\mathfrak{g}}$ is the quotient field of $Y(\mathfrak{g})$.  It is the field generated by $a_{j_1},\ldots, a_{j_r}$, which are algebraically independent over $k$.  In particular, $r = i(\mathfrak{g})$.
\item[(4)] $Y(\mathfrak{g}) \subset k[a_{j_1},\ldots, a_{j_r}, a^{-1}]$ for some nonzero $a \in k[a_{j_1},\ldots, a_{j_r}]$ 
\end{itemize}

{\bf Theorem 45.}\\
Any nilpotent Lie algebra $\mathfrak{g}$ with $i(\mathfrak{g}) \leq 2$ is coregular.\\
\ \\
{\bf Proof.} We observe that
$$1 \leq \dim Z(\mathfrak{g}) \leq i(\mathfrak{g}) \leq 2$$
If $\dim Z(\mathfrak{g}) = i(\mathfrak{g})$ then $\mathfrak{g}$ is square integrable and therefore $Y(\mathfrak{g}) = S(Z(\mathfrak{g}))$ which is polynomial by Remark 38
(or section 2.8).  So it suffices to deal with the case where $\dim Z(\mathfrak{g}) = 1$ and $i(\mathfrak{g}) = 2$. (hence $\dim \mathfrak{g}$ is even)\\
Choose a sequence of ideals
$$0 = \mathfrak{g}_0 \subset \mathfrak{g}_1 \subset \ldots \subset \mathfrak{g}_n = \mathfrak{g}$$
with the same properties as in Theorem 44.  In particular, $x_1,\ldots, x_n$ is a basis of $\mathfrak{g}$ with $x_i \in \mathfrak{g}_i\backslash \mathfrak{g}_{i-1}$, $i = 1,\ldots, n$.  Note that $\mathfrak{g}_1 = \langle x_1\rangle = Z(\mathfrak{g})$.  So,
$$x_1 \in S(\mathfrak{g}_1) \cap\ Y(\mathfrak{g}) = k[x_1],\ \mbox{but}\ x_1 \notin S(\mathfrak{g}_0) \cap\ Y(\mathfrak{g}) = k$$
Since $i(\mathfrak{g}) = 2$ there is only one more $j > 1$ such that
$$S(\mathfrak{g}_{j-1}) \cap  Y(\mathfrak{g}) \subs\limits^{}_{\neq} S(\mathfrak{g}_j) \cap\ Y(\mathfrak{g})$$
and there is a nonzero.
$$b \in S(\mathfrak{g}_{j-1}) \cap \ Y(\mathfrak{g}) = S(\mathfrak{g}_1)\cap  Y(\mathfrak{g}) = k[x_1]\ \ \ \mbox{and}\ \ \  c \in S(\mathfrak{g}_{j-1})$$
such that
$$a = bx_j + c \in S(\mathfrak{g}_j) \cap\ Y(\mathfrak{g})$$
by Theorem 44.  We may assume that $a$ is of smallest degree with these properties.  Clearly the homogeneous components of $a$ also belong to $S(\mathfrak{g}_j) \cap\ Y(\mathfrak{g})$.  Among them
we let $f$ be a homogeneous component not contained in $S(\mathfrak{g}_{j-1}) \cap\ Y(\mathfrak{g}) = k[x_1]$.  Note that $\deg f = \deg a$.  Then up to a nonzero scalar multiplier:\\
$f = x_1^m x_j + u$ for some $m$ and homogeneous $u \in S(\mathfrak{g}_{j-1})$ of degree $m + 1$.\\
Clearly $u \neq 0$ (otherwise $x_1^m x_j \in Y(\mathfrak{g})$ and hence $x_j \in Z(\mathfrak{g}) = \mathfrak{g}_1$, contradiction)\\
Also, $m \geq 1$ (if $m=0$ then $u \in \mathfrak{g}_{j-1}$ and $f = x_j + u \in Z(\mathfrak{g}) = \mathfrak{g}_1 = \langle x_1\rangle$, contradiction).\\
Moreover, $u$ is not divisible by $x_1$ (otherwise $f/x_1$, which belongs to $S(\mathfrak{g}_j) \cap\ Y(\mathfrak{g})$ but not to $S(\mathfrak{g}_{j-1})\cap \ Y(\mathfrak{g}) = k[x_1]$, would be of a smaller degree than $f$, contradiction).\\
\ \\
{\bf Claim.}  Suppose we have an element $f = x_1^m x_j + u \in Y(\mathfrak{g})$ with $m \geq 1$, $j > 1$ and nonzero homogeneous $u \in S(\mathfrak{g}_{j-1})$ of degree $m+1$ 
and not divisible by $x_1$.  Then $Y(\mathfrak{g}) = k[x_1, f]$, which is polynomial.\\
\ \\
First we see that $f \in S(\mathfrak{g}_j) \cap\ Y(\mathfrak{g})$, $x_1^m \in S(\mathfrak{g}_{j-1}) \cap \ Y(\mathfrak{g}) = k[x_1]$ and $u \in S(\mathfrak{g}_{j-1})$. By Theorem 44, $x_1$ and $f$ are algebraically independent
over $k$ and
$$Y(\mathfrak{g}) \subset k[x_1, f, x_1^{-m}] \subset k[x_1, f, x_1^{-1}]$$
So, it suffices to show that $Y(\mathfrak{g}) = k[x_1, f]$ ($\ast$). For this we need the following lemmas:\\
\\
{\bf Lemma A:} Let  $P \in k[X]$ be a polynomial.  If $x_1$ divides $P(f)$ then $P = 0$ and so $P(f) = 0$.\\
\ \\
{\bf Proof. } Let $I$ be the ideal of $S(L)$ generated by $x_1$.  We identify the quotient algebra by $k [x_2,\ldots, x_n]$.  Denote by $f_1$ and $u_1$ the canonical images of $f$ and $u$.  By assumption we obtain $P(u_1) = P(f_1) = 0$.\\
But $u_1 \in k[x_2,\ldots, x_n]$ is nonzero (as $x_1$ does not divide $u$) of degree $m+1$.  Consequently, $P = 0$ and $P(f) = 0$. \hfill $\square$\\
\ \\
{\bf Lemma B.}   Let $q \in S(\mathfrak{g})$ be such that $x_1 q \in k[x_1, f]$.  Then also $q \in k[x_1, f]$.\\
\ \\
{\bf Proof.} By assumption there are $h_i \in k[f]$ such that, for some $r$
$$x_1 q = x_1^r h_r + \ldots + x_1^2 h_2 + x_1 h_1 + h_0$$
Clearly, $x_1$ divides $h_0$, which implies that $h_0 = 0$ by Lemma A.  Consequently,
$$q = x_1^{r-1} h_r + \ldots + x_1 h_2 + h_1 \in k[x_1, f]$$
\hfill $\square$\\
\ \\
We can now show ($\ast$) as follows:\\
Take $q \in Y(\mathfrak{g}) \subset k [x_1, f, x_1^{-1}]$, i.e. $x_1^t q \in k[x_1,f]$ for some $t$.\\
By applying Lemma B $t$ times we may conclude that $q \in k[x_1, f]$.  Therefore $Y(\mathfrak{g}) \subset k[x_1, f]$.  The other inclusion is obvious.\\
Finally, we provide a formula for $\deg f$.  First, $c(\mathfrak{g}) = \ds\frac{1}{2} (\dim \mathfrak{g} + i(\mathfrak{g})) = \ds\frac{n}{2} + 1$. Then, by Theorem 18:
$$\deg x_1 + \deg f = c(\mathfrak{g}) - \deg p_\mathfrak{g} = \ds\frac{n}{2} + 1 - \deg p_\mathfrak{g}$$
Consequently, $\deg f = \ds\frac{n}{2} - \deg p_\mathfrak{g}$ (in particular, $\deg f = \ds\frac{n}{2}$ if $\mathfrak{g}$ is nonsingular)\hfill $\square$\\
\ \\
{\bf Remark 46.}  Yakimova has informed us that Theorem 45 can also be derived from a result of Michel Brion on linear actions of unipotent groups [Br].\\
\ \\
{\bf Remark 47.} Theorem 45 fails if
\begin{itemize}
\item[(1)] the condition on $i(\mathfrak{g})$ is replaced by $i(\mathfrak{g}) = 3$. \\
Indeed, the standard filiform Lie algebra $\mathfrak{g}_{5,5}$ has index 3, but is not coregular [O7, p. 1304].  This was already known by Dixmier [D2, Proposition 2].
\item[(2)] nilpotent is replaced by solvable.\\
Indeed, Example 23 is solvable of index 2, but it is not coregular.
\end{itemize}
{\bf Question 48.} Let $\mathfrak{g}$ be nilpotent with $i(\mathfrak{g}) \leq \dim Z(\mathfrak{g}) + 1$.  Does this imply that $\mathfrak{g}$ is coregular ?\\
This is true for all indecomposable nilpotent Lie algebras of dimension at most seven [O7, O8].\\
\ \\
We will now examine the coregularity of the major types of filiform Lie algebras, presented in [GK1, p.41].\\
\ \\
{\bf Definition 49.}\\
Consider the descending central series of $\mathfrak{g}$
$$C'(\mathfrak{g}) = \mathfrak{g}, C^2(\mathfrak{g}) = [\mathfrak{g},\mathfrak{g}],\ldots, C^i (\mathfrak{g}) = [\mathfrak{g}, C^{i-1}(\mathfrak{g})],\ldots$$
An $n$-dimensional Lie algebra $\mathfrak{g}$ is called filiform if $\dim C^i (\mathfrak{g}) = n-i$, $i = 2,\ldots, n$.\\
In particular, $C^n(\mathfrak{g}) = 0$ (and thus $\mathfrak{g}$ is nilpotent) and $Z(\mathfrak{g}) = C^{n-1}(\mathfrak{g})$ is $1$-dimensional.\\
\ \\
Combining Proposition 42 with Theorem 45 yields:\\
\ \\
{\bf Proposition 50.} Let $\mathfrak{g}$ be an 8-dimensional filiform Lie algebra.  Then
\begin{center}
$\mathfrak{g}$ is coregular\ \ \ $\Leftrightarrow \ \ \ i(\mathfrak{g}) = 2$
\end{center}

{\bf Theorem 51.} Let $\mathfrak{g}$ be an $n$-dimensional filiform Lie algebra.  Then
\begin{itemize}
\item[(1)] If $\mathfrak{g}$ is of type $Q_n$ or $W_n$ then $i(\mathfrak{g}) \leq 2$, so $\mathfrak{g}$ is coregular.
\item[(2)] If $\mathfrak{g}$ is of type $L_n$ or $R_n$ then
\begin{center}
$\mathfrak{g}$ is coregular \ \ \ $\Leftrightarrow \ \ \ i(\mathfrak{g}) \leq 2$
\end{center}
\end{itemize}

{\bf Proof.}
\begin{itemize}
\item[(1)] a) Suppose $\mathfrak{g}$ is of type $Q_n$.\\
Basis of $\mathfrak{g}:\ x_1,\ldots, x_n$, $n = 2q$.\\
Nonzero brackets: \ $[x_1,x_i] = x_{i+1}$, $i = 2,\ldots, n-2$ and $[x_j, x_{n-j+1}] = (-1)^{j+1} x_n$, $j = 2,\ldots, q$.\\
We observe that
$$i(\mathfrak{g}) = 2, \ \ \ Z(\mathfrak{g}) = \langle x_n\rangle,\ \ \ p_\mathfrak{g} = x_n^{q-2}$$
So, $\mathfrak{g}$ is coregular by Theorem 45.  Next put
$$f = 2x_1x_n + (-1)^{q+1} x_{q+1}^2 + 2\sum\limits_{i=3}^q (-1)^i x_i x_{n-i+2} \in Y(\mathfrak{g})$$
which satisfies the conditions of the claim in the proof of Theorem 45.  Therefore $Y(\mathfrak{g}) = k[x_n,f]$.\\
\ \\
As an alternative solution we can use the Frobenius method:\\
Consider the torus $T = \langle t_1, t_2\rangle \subset \Der \ \mathfrak{g}$, where 
$$t_1 = \mbox{diag} (0, 1,1,\ldots, 1,2),\ \ \ t_2 = \mbox{diag} (1,1,2,3,\ldots, n-2, n-1)$$
see [R, p.4].  Then the semi-direct product $L = T \oplus \mathfrak{g}$ is Frobenius.  Indeed,
$$\Delta (L) = x_n^{n-2} f^2 \neq 0. \ \ \mbox{Hence}\ \ p_L = x_n^{q-1}f$$
By (4) of Theorem 30:
$$Y(\mathfrak{g}) = Sy(L) = k[x_n,f]$$
since $\mathfrak{g}$ is nilpotent (and hence has no proper semi-invariants) and so $\mathfrak{g} = L_\Lambda$.\\
\ \\
b) Suppose $\mathfrak{g}$ is of type $W_n$\\
Basis: \ \ \ $x_1, \ldots, x_n$\\
Nonzero brackets: \ \  $[x_i, x_j] = (j-i) x_{i+j}$, $i < j$ and $i+j \leq n$
\begin{itemize}
\item[b1)] $n = 2q+1$.  Then $Z(\mathfrak{g}) = \langle x_n\rangle$ and $\dim Z(\mathfrak{g}) = 1 = i(\mathfrak{g})$, i.e. $\mathfrak{g}$ is square integrable.  Consequently, $Y(\mathfrak{g}) = k[x_n]$ by Remark 38.
\item[b2)] $n = 2q$.  Then $i(\mathfrak{g}) = 2$ and thus $\mathfrak{g}$ is coregular by Theorem 45.\\
For example, if $ n = 8$ then $Y(\mathfrak{g}) = k[x_8, f]$ where
$$ f = 64x_4 x_8^3 - 16x_6^2x_8^2 - 32 x_5 x_7 x_8^2 + 24 x_6 x_7^2 x_8 - 5x_7^4$$
\end{itemize}
\item[(2)] i) Suppose $\mathfrak{g}$ is of type $L_n$, the standard filiform Lie algebra,\\
Basis : $x_1,\ldots, x_n$, $n \geq 3$\\
Nonzero brackets: $[x_1, x_i] = x_{i+1}$, $i=2,\ldots, n-1$.\\
Clearly, $i(\mathfrak{g}) = n-2$ and $\mathfrak{h} = \langle x_2, \ldots, x_n\rangle$ is an abelian ideal of codimension one of $\mathfrak{g}$.  Also, $[\mathfrak{g},\mathfrak{g}] = \langle 
x_3,\ldots, x_n\rangle$ and $\dim[\mathfrak{g},\mathfrak{g}] = n-2$.\\
By Theorem 20:
\begin{center}
$\mathfrak{g}$ is coregular \ \ \ $\Leftrightarrow \ \ \ \dim[\mathfrak{g},\mathfrak{g}] \leq 2\ \ \ \Leftrightarrow \ \ \ n \leq 4 \ \ \ \Leftrightarrow\ \ \ i(\mathfrak{g}) \leq 2$
\end{center}
See also [OV, Example 1.7] and [O8, Example 27].\\[0ex]
[The Lie algebras 1,2,8,25,159 of [O7, O8] are of this type].
\item[ii)]  Suppose $\mathfrak{g}$ is of type $R_n$.\\
Basis: $x_1, \ldots, x_n$, $n \geq 5$.\\
Nonzero brackets: $[x_1, x_i] = x_{i+1}$, $i = 2,\ldots, n-1$; $[x_2, x_j] = x_{j+2}$, $j = 3,\ldots, n-2$.\\
Since $i(\mathfrak{g}) = n-4$, it suffices to show that
$$\mathfrak{g}\ \mbox{is coregular}\ \ \ \Leftrightarrow \ \ \ n \leq 6$$
First suppose $\mathfrak{g}$ is coregular.  Note that $c(\mathfrak{g}) = \ds\frac{1}{2} (n + n -4) = n-2$.  Then $\mathfrak{h} = \langle x_3, x_4,\ldots, x_n\rangle$ is an 
$(n-2)$-dimensional abelian ideal of $\mathfrak{g}$ and so is a CP of $\mathfrak{g}$.  Next, it is not difficult to see that
$$\mathfrak{g}_{\sing}^\ast = \{f \in \mathfrak{g}^\ast \mid f(x_5) = \ldots = f(x_n) = 0\}$$
Hence, $\codim \ \mathfrak{g}_{\sing}^{\ast} = n-4$.  By Theorem 16 $\codim \ \mathfrak{g}_{\sing}^\ast \leq 2$, i.e. $n \leq 6$.\\
Conversely, if $n \leq 6$ then $i(\mathfrak{g}) \leq 2$ and so $\mathfrak{g}$ is coregular by Theorem 45. [The Lie algebras 6,27,151 of [O7, O8] are of this type] \hfill $\square$
\end{itemize}

{\bf Remark.} There are coregular filiform Lie algebras of index larger than 2.  For instance the Lie algebra 106 of [O8, p.104].\\
\ \\
{\bf Theorem 52.} Let $\mathfrak{g}$ be a quadratic Lie algebra.  Then $\mathfrak{g}$ is coregular if one of the following conditions is satisfied:
\begin{itemize}
\item[(i)] $[\mathfrak{g},\mathfrak{g}]\neq \mathfrak{g}$ and $i(\mathfrak{g}) = 2$
\item[(ii)] $\mathfrak{g}$ is nilpotent and $i(\mathfrak{g}) = 3$
\end{itemize}
\ \\
{\bf Proof.} 
\begin{itemize}
\item[(i)] Since $i(\mathfrak{g}) = 2$ we have that $ n = \dim \mathfrak{g}$ is even and $c(\mathfrak{g}) = \ds\frac{n+2}{2} = \ds\frac{n}{2} + 1$.  
We may assume that $n \geq 4$. $\mathfrak{g}$ being quadratic, admits a nondegenerate, symmetric, invariant bilinear form $b$ 
(such a Lie algebra is sometimes called regular quadratic) [FS].  It is easy to verify that w.r.t. $b$ we obtain that
$$Z(\mathfrak{g}) = [\mathfrak{g},\mathfrak{g}]^\perp \neq 0\ \mbox{since}\ \ [\mathfrak{g},\mathfrak{g}] \neq \mathfrak{g}$$
$\mathfrak{g}$ is a fortiori quasi quadratic (2.8), i.e.
$$\mathfrak{g} = F(\mathfrak{g}) = \sum\limits_{\xi \in \mathfrak{g}_{\reg}^\ast} \mathfrak{g}(\xi)$$
In particular, $\mathfrak{g}(\xi) \neq \mathfrak{g}(\eta)$ for some $\xi, \eta \in \mathfrak{g}_{\reg}^\ast$, which we extend to algebra endomorphisms of $S(\mathfrak{g})$.  Both $\mathfrak{g}(\xi)$ and $\mathfrak{g}(\eta)$ contain $Z(\mathfrak{g})$ and are of dimension $i(\mathfrak{g}) = 2$.  It follows that $\dim Z(\mathfrak{g}) = 1$.  Hence, we can find a basis $x_1, x_2, \ldots, x_{n-1}, x_n$ of $\mathfrak{g}$ such that
$$\mathfrak{g}(\xi) = \langle x_{n-1}, x_n\rangle, \ \ \ \mathfrak{g}(\eta) = \langle x_1,x_n\rangle \ \ \ \mbox{and}\ \ \ Z(\mathfrak{g}) =  \langle x_n\rangle$$
By 2.1 the rank $r$ of the structure matrix $M = ([x_i, x_j])$ is given by
$$r = \dim \mathfrak{g} - i(\mathfrak{g}) = n-2$$
Clearly, $\langle x_1, \ldots, x_{n-2}\rangle \oplus \mathfrak{g}(\xi) = \mathfrak{g}$.\\
Next, we consider the $r \times r$ submatrix $A = ([x_i,x_j])_{1 \leq i,j \leq n-2}$ with Pfaffian $p \in S(\mathfrak{g})$.  Then $\xi(p) \neq 0$, indeed
$$\xi(p)^2 = \xi(p^2) = \xi(\det A) = \det (\xi([x_i,x_j]) \neq 0$$
since $\xi$ is regular.  Similarly, we  observe that
$$\langle x_2 ,\ldots, x_{n-1}\rangle \oplus \mathfrak{g}(\eta) = \mathfrak{g}$$
We put $B = ([x_i, x_j])_{2 \leq i,j\leq n-1}$ with Pfaffian $q \in S(\mathfrak{g})$.  As before we get $\eta (q) \neq 0$.  On the other hand, $\eta(p) = 0$ because
$$\eta(p)^2 = \eta (p^2) = \eta (\det A) = \det (\eta([x_i,x_j])_{1 \leq i,j\leq n-2}) = 0$$
since the first row of the matrix is
$$(\eta([x_1,x_1]),\eta([x_1, x_2]),\ldots, \eta([x_1, x_{n-2}]))$$
which is zero because $x_1 \in \mathfrak{g}(\eta)$.\\
Consequently, $p$ and $q$ are principal $r \times r$ Pfaffians of the structure matrix $M$ of $\mathfrak{g}$ of the same degree, namely $\ds\frac{n-2}{2} = \ds\frac{n}{2} - 1$.  By the above $p$ is not a scalar multiple of $q$.  The fundamental semi-invariant $p_\mathfrak{g}$, being the GCD of all principal $r \times r$ Pfaffians (Definition 5), divides both $p$ and $q$.  Therefore,
$$\deg p_\mathfrak{g} \leq (\ds\frac{n}{2} - 1) - 1 = \ds\frac{n}{2}- 2$$
Next, let $y_1,\ldots, y_n$ be the dual basis of $x_1, \ldots, x_n$ w.r.t. $b$, i.e. $b(x_i,y_j) = \delta_{ij}$ for all $i, j : 1,\ldots, n$.  Then,
$$ f= x_1y_1 + \ldots + x_n y_n \in Y(\mathfrak{g})$$
is the well known Casimir element of $S(\mathfrak{g})$. Finally, $x_n$ and $f$ are algebraically independent, homogeneous elements of $Y(\mathfrak{g})$ such that
$$\deg x_n + \deg f = 3 = (\ds\frac{n}{2} + 1) - (\ds\frac{n}{2} - 2) \leq c(\mathfrak{g}) - \deg p_\mathfrak{g}$$
By Theorem 29 equality holds (in particular $\deg p_\mathfrak{g} = \ds\frac{n}{2} - 2$) and \\
$Y(\mathfrak{g}) = k[x_n, f]$. 
\item[(ii)] Since $i(\mathfrak{g}) = 3$ we have that $n = \dim\ \mathfrak{g}$ is odd and $c(\mathfrak{g}) = \ds\frac{1}{2}(n+3)$.  We may assume that $n \geq 5$ (otherwise $\mathfrak{g}$ is abelian).
$\mathfrak{g}$ being quadratic, admits a nondegenerate, symmetric, invariant bilinear form $b$.\\
$\mathfrak{g}$ is a fortiori quasi quadratic, i.e.
\begin{eqnarray*}
\mathfrak{g} = F(\mathfrak{g}) = \sum\limits_{\xi \in \mathfrak{g}_{reg}^\ast}  \mathfrak{g}(\xi)
\end{eqnarray*}
In particular, $\mathfrak{g}(\xi) \neq \mathfrak{g}(\eta)$ for some $\xi, \eta \in \mathfrak{g}_{\reg}^\ast$ which we extend to algebra endomorphisms of $S(\mathfrak{g})$.  Both $\mathfrak{g}(\xi)$ and $\mathfrak{g}(\eta)$ contain $Z(\mathfrak{g})$ and are of dimension $i(\mathfrak{g}) = 3$.\\
On the other hand, $\dim Z(\mathfrak{g}) \geq 2$ because $\mathfrak{g}$ is nilpotent and quasi quadratic [O5, Corollary 3.6].  Therefore $\dim Z(\mathfrak{g}) = 2$.  Hence there exists a basis
$x_1, x_2,\ldots, x_{n-2}, x_{n-1}, x_n$ of $\mathfrak{g}$ such that $\mathfrak{g}(\xi) = \langle x_{n-2}, x_{n-1}, x_n\rangle$, $\mathfrak{g}(\eta) = \langle x_1, x_{n-1}, x_n\rangle$ and 
$Z(\mathfrak{g}) = \langle x_{n-1}, x_n\rangle$.
The rank $r$ of the structure matrix $M = ([x_i,x_j])$ is given by
$$r = \dim\ \mathfrak{g} - i(\mathfrak{g}) = n-3$$
Similar to the proof of (i) we can find principal $r \times r$ Pfaffians $p$ and $q$ of degree $\ds\frac{1}{2}(n-3)$ such that one is not a scalar multiple of the other.\\
Since the fundamental semi-invariant $p_{\mathfrak{g}}$ divides both $p$ and $q$, we get
$$\deg p_{\mathfrak{g}} \leq \ds\frac{1}{2} (n-3) - 1 = \ds\frac{1}{2}(n-5)$$
Next, let $f \in S(\mathfrak{g})$ be the Casimir element w.r.t. $b$. Then $x_{n-1}, x_n, f$ are algebraically independent, homogeneous elements of $Y(\mathfrak{g})$ such that
$$\deg x_{n-1} + \deg x_n + \deg f = 4 = \ds\frac{1}{2}(n+3) - \ds\frac{1}{2} (n-5) \leq c(\mathfrak{g}) - \deg p_{\mathfrak{g}}$$
By Theorem 29 equality holds (in particular $\deg p_{\mathfrak{g}} = \ds\frac{1}{2} (n-5)$) and $Y(\mathfrak{g}) = k[x_{n-1}, x_n,f]$. \hfill $\square$\\
\ \\
\end{itemize}

{\bf 5. Polynomiality for nonsolvable Lie algebras of dimension at most eight}\\
Because of Example 23 we now restrict ourselves to the nonsolvable case.\\
\ \\
{\bf Theorem 53.}  Let $L$ be a nonsolvable, indecomposable Lie algebra with $\dim L \leq 8$.  Then
\begin{itemize}
\item[1)] $Y(L)$ and $Sy(L)$ are polynomial algebras over $k$ (and hence so are $Z(U(L))$ and $Sz(U(L))$)
\item[2)] $R(L)^L$ (and hence also $Z(D(L))$) is rational over $k$.
\end{itemize}

{\bf Proof.}  This will proceed case by case using the classification provided to us by B. Komrakov.  See also [Tu]. In [AOV2] the algebraic ones among them were shown to satisfy the following well known\\
{\bf Gelfand-Kirillov conjecture [GK]}. Let $L$ be an algebraic Lie algebra over $k$.  Then $D(L)$ is isomorphic to a Weyl skew field $D_n(F)$ over a rational extension $F$ of $k$.  In particular, $Z(D(L))$, which is isomorphic to $F$, is also rational over $k$.\\
\ \\
Over the years positive, but also some negative, answers have been obtained [BGR, J1, Mc, N, AOV1, AOV2, O6, Pr].  See also Appendix.\\
As we will use some results of [AOV2], we will employ the same notation.
\\
Before we list in 5.1 for each case the results of the verification of the theorem, we would like to present in detail some typical examples in order to exhibit the various procedures
used in the proof.\\
\ \\
{\bf Example 54.} Let $L$ be the semi-direct product $L_{6,3} = sl(2,k)\oplus H$ of $sl(2,k)$ with the 3-dimensional Heisenberg Lie algebra $H$ with basis 
$h, x, y, e_0, e_1, e_2$
and nonzero brackets:\\
$[h,x] = 2x$, $[h,y] = -2y$, $[x,y] = h$, $[h, e_0] = e_0$, $[h,e_1] = -e_1$, $[x,e_1] = e_0$, $[y, e_0] = e_1$, $[e_0, e_1] = e_2$\\
Clearly, $i(L) = 2$ and $c(L) = \ds\frac{1}{2} (6+2) = 4$.  One verifies that $p_L = 1$ and also that $L$ is quasi quadratic (i.e. $F(L) = L$), so there are no CP's.  Next, $Y(L)$
contains the following homogeneous, algebraically independent elements:
$$e_2\ \ \ \mbox{and}\ \ \ f= e_2(h^2 + 4xy) + 2(e_0e_1h + e_1^2x - e_0^2y)$$
Because $\deg e_2 + \deg f = 4 = c(L) - \deg p_L$ we may conclude that $Y(L) = k[e_2, f]$ by Theorem 29.  Since $[L,L] = L$, $L$ is algebraic, without proper semi-invariants.
Therefore
$$Sy(L) = Y(L)\ \ \ \mbox{and}\ \ \ R(L)^L = k(e_2, f)$$
Finally, $M = k[e_1, e_2, e_0^2 - 2e_2x, f]$ is a polynomial, complete, Poisson commutative subalgebra of $S(L)$.\\
\ \\
{\bf Example 55.} Let $L$ be the 7-dimensional algebraic Lie algebra $L_{7,9}$ with basis
$h, x, y, e_0, e_1, e_2, e_3$ and with nonzero brackets:\\
$[h,x] = 2x$, $[h,y] = -2y$, $[x,y] = h$, $[h, e_0] = e_0$, $[h,e_1] = -e_1$, $[x,e_1] = e_0$, $[y, e_0] = e_1$, $[e_0, e_1] = e_2$, $[e_3, e_0] = e_0$, $[e_3, 
e_1] = e_1$, $[e_3, e_2] = 2e_2$.\\
Note that $L$ is the semi-direct product of $L_{6,3}$ with $\ad\ e_3$, which is semi-simple.  One verifies that $i(L) = 1$, $c(L) = 4$ and $p_L = 1$.  Also
$$F(L) = \langle h, x, y, e_0, e_1, e_2\rangle = [L,L] = L_{6,3}$$
Since this is not commutative there are no CP's.  The above implies that $L_\Lambda = [L,L]$.  By (5) of Theorem 4 we observe that
$$Sy(L) = Y(L_\Lambda) = Y(L_{6,3}) = k[e_2, f]$$
where $f = e_2(h^2 + 4xy) + 2(e_0e_1 h + e_1^2x - e_0^2y)$.\\
Furthermore, $e_2$ and $f$ are irreducible semi-invariants with the same weight $\lambda \in L^{\ast}$ for which $\lambda(e_3) = 2$ and $\lambda([L,L]) = 0$.\\
Consequently, $Y(L) = k$ and $R(L)^L = k(e_2^{-1}f)$.\\
Finally, we know from the previous example that $M = k[e_1, e_2, e_0^2 - 2e_2x, f]$ is a polynomial, complete, Poisson commutative subalgebra of $S(L_{6,3})$.  
Hence $\trdeg_kM = c(L_{6,3}) = c(L_\Lambda) = c(L)$,
the latter by (4) of Theorem 4.  Therefore, $M$ is also complete in $S(L)$.\\
\ \\
{\bf Example 56.} Let $\mathfrak{g}$ be the 8-dimensional algebraic Lie algebra $L_{8,17}$ with basis $h, x, y, e_0, e_1, e_2, e_3, e_4$ and nonzero brackets:\\
$[h,x] = 2x$, $[h,y] = -2y$, $[x,y] = h$, $[h, e_0] = e_0$, $[h,e_1] = -e_1$, $[h,e_2] = e_2$, $[h,e_3] = -e_3$, $[x,e_1] = e_0$, $[x,e_3] = e_2$, $[y, e_0] = e_1$, 
$[y,e_2] = e_3$, $[e_2, e_4] = e_0$, $[e_3, e_4] = e_1$.\\
One verifies that $i(\mathfrak{g}) = 2$, $c(\mathfrak{g}) = 5$, $p_\mathfrak{g} = 1$.  Also $\mathfrak{g}$ is quasi-quadratic (i.e. $F(\mathfrak{g}) = \mathfrak{g})$ and so
it has no proper semi-invariants (2.8).  Moreover, there are no CP's as $F(\mathfrak{g})$ is not commutative (2.8).\\
Next, consider the semi-direct product $L = T \oplus \mathfrak{g}$, where $T = \langle t_1, t_2\rangle \subset \Der\ \mathfrak{g}$ with
$$t_1 = \mbox{diag} (0,0,0,1,1,0,0,1),\ \ \ t_2 = \mbox{diag} (0,0,0,1,1,1,1,0)$$
Then
$$\Delta(L) = 4(e_1e_2 - e_0e_3)^2 (e_0e_1h + e_1^2x - e_0^2y + e_1e_2e_4 - e_0 e_3 e_4)^2 \neq 0$$
Hence, $L$ is a 10-dimensional Frobenius Lie algebra.  By (4) of Theorem 30 we may conclude that
$$Y(\mathfrak{g}) = k[f_1,f_2] = Sy(\mathfrak{g})\ \ \ \mbox{and}\ \ \ R(\mathfrak{g})^\mathfrak{g}= k(f_1, f_2)$$
where $f_1 = e_1e_2 - e_0e_3$ and $f_2 = e_0 e_1h + e_1^2x - e_0^2y + (e_1e_2 - e_0e_3)e_4$.\\
Finally, $M = k[e_0, e_1, e_2, e_3, f_2]$ is a polynomial, complete Poisson commutative subalgebra of $S(\mathfrak{g})$.\\
\ \\
{\bf Remark 57.} The same Frobenius method can be used to show the theorem for $L_{8,2}$ [O7, p.1302].  The theorem also holds for $L_{6,4}$, $L_{8,19}$, $L_{8,20} (\alpha \neq -1)$, $L_{8,28}$ (since they
 are Frobenius) as well as for their canonical truncations $L_5$, $L_{7,1}$, $L_{7,2}$ (apply Theorem 10).\\
 \ \\
{\bf Example 58.} Let $L$ be the 8-dimensional non algebraic Lie algebra $L_{8,25}$ with basis $h, x, y, e_0, e_1, e_2, e_3, e_4$ and nonzero brackets:\\
$[h,x] = 2x$, $[h,y] = -2y$, $[x,y] = h$, $[h, e_0] = e_0$, $[h,e_1] = -e_1$, $[x,e_1] = e_0$, $[y, e_0] = e_1$, $[e_4, e_0] = e_0$, $[e_4, e_1] = e_1$, $[e_4, e_2] = -e_3$.\\
One verifies that $i(L) = 2$, $c(L) = 5$, $p_L = 1$, $Z(L) = \langle e_3\rangle$ and
$$F(L) = \langle h, x, y, e_0, e_1, e_2, e_3\rangle = L_\Lambda$$
Next, we put $E_1 = \ad\ h$, $E_2 = \ad\ x$, $E_3 = \ad\ y$, $E_4 = \ad\ e_0$, $E_5 = \ad\ e_1$, $E_6 = \ad\ e_2$, $E_7(e_0) = e_0$, $E_7 (e_1) = e_1$ and zero on others, $E_8(e_2) = -e_3$ and zero
on others.  Then
$$\ad\ e_4 = E_7 + E_8$$
is the decomposition of $\ad\ e_4$ into its semi-simple and nilpotent components.  It follows that $E_1, E_2,\ldots, E_8$ is a basis for the algebraic hull $H$ of $\ad\ L$.\\
Next, we rename the basis
$$h,x,y,e_0, e_1,e_2,e_3, e_4\ \ \ \mbox{by}\ \ \ x_1, x_2, x_3,\ldots, x_8$$
Then we see that $\rank(E_ix_j) = 7$ and so
$$j(L) = \dim L - \rank(E_ix_j) = 1 = \dim Z(L)$$
So, we can apply Theorem 37. First one verifies that $p'_L = e_3 (e_0e_1h + e_1^2x - e_0^2y)$.  
\\
Hence, $f = e_0e_1h + e_1^2x - e_0^2y$ is the only proper irreducible semi-invariant and we may conclude that
$$ Y(L) = k[e_3],\ \ \ R(L)^L = k(e_3),\ \ \ Sy(L) = Y(L)[f] = k[e_3, f]$$
Finally, $M = k[e_0, e_1, e_2, e_3, f]$ is a polynomial, complete, Poisson commutative subalgebra of $S(L)$.\\
\ \\
The result on $Sy(L)$ can also be seen as follows.  Clearly, $L = L_\Lambda \oplus ke_4$ while $L_\Lambda = L_5 \times \langle e_2, e_3\rangle$ (direct product).\\
Now $Y(L_5) = k[f]$ [O7, p.1301].\\
By (5) of Theorem 4:
$$Sy(L) \subset Y(L_\Lambda) = k[e_2, e_3, f]$$

On the other hand, $e_3$ is an invariant and $f$ is a semi-invariant for $L$, indeed
$$\{x,f\} = \lambda(x)f\ \mbox{for all}\ x \in L$$
where $\lambda \in L^{\ast}$, $\lambda(e_4) = 2$ and zero on others.  Consequently,
$$k[e_3, f] \subset Sy(L) \subset k[e_2, e_2, f]$$
Now take any semi-invariant $g \in Sy(L)$.  Then already $g \in k[e_2, e_3, f]$.  As $g$ is a semi-invariant for $\ad\ L$ it is also one under 
the action of $H$ [C, p.208].  In particular,
$$E_8(g) = ag\ \mbox{for a suitable}\ a\in k$$
But $E_8$ is nilpotent and so $a=0$.\\
Next, we consider
$$-e_3 \ds\frac{\partial g}{\partial e_2} = -x_7 \ds\frac{\partial g}{\partial x_6} = \sum\limits_{j=1}^8 E_8 (x_j) \ds\frac{\partial g}{\partial x_j} = E_8(g) = 0$$
which implies that $\ds\frac{\partial g}{\partial e_2} = 0$ and so $g \in k[e_3, f]$.  Therefore $Sy(L) = k[e_3, f]$.\\
\ \\
{\bf 5.1. List of indecomposable nonsolvable Lie algebras of dimension $\vet{\leq 8}$}\\
The main purpose is to show that for each member $L$ of the list $Y(L)$, $Sy(L)$ and $R(L)^L$ satisfy the requirements of Theorem 53 by giving their explicit description.  In particular,
 $L$ is coregular.  It will turn out that $i(L) \leq 2$, which is hardly surprising in view of Proposition 42. In addition we will provide the Frobenius semi-radical $F = F(L)$ and if it exists a CP-ideal (CPI).\\
 $M$ will be a polynomial, complete, Poisson commutative subalgebra of $S(L)$. Other abbreviations are: $i = i(L)$, $c = c(L)$, $p = p_L$, $Y = Y(L)$, $Sy = Sy(L)$, $R^I = R(L)^L$.  Furthermore,
 $h, x, y$ will be the standard basis of $sl(2, k)$.  $W_n$ will be its $(n+1)$-dimensional  irreducible module with standard basis $e_0, e_1,\ldots, e_n$.  In particular, 
 $$h\cdot e_i = (n-2i) e_i,\ x\cdot e_i = (n-i+1) e_{i-1},\ y\cdot e_i = (i+1) e_{i+1}$$
 for all $i$ and $e_{-1} = e_{n+1} = 0$.\\
 \ \\
 {\bf I. $\vet{L}$ is algebraic}\\
 For this the possible parameters need to be rational numbers, but we will briefly indicate what happens if they are not.\\
 \ \\
 {\bf I.0. $\vet{\dim L = 3}$}
 \begin{itemize}
 \item[0.] $sl (2,k)$ (simple and hence quadratic)\\
 Basis: $h, x, y$\\
 $[h,x] = 2x$, $[h,y] = -2y$, $[x,y] = h$.\\
 $i=1$, $c=2$, $p=1$, $F = sl (2,k)$, no CP's,\\
 $Y = k[f] = Sy$, where $f = h^2 + 4xy$\\
 $R^I = k(f)$, $M = k[h,f]$.
 \end{itemize}
 {\bf I.1.  $\vet{\dim L = 5}$}
 \begin{itemize}
 \item[1.] $L_5 = sl(2,k) \oplus W_1$ (quasi quadratic)\\
 Basis: $h, x, y,, e_0, e_1$\\
 $[h,x] = 2x$, $[h,y] = -2y$, $[x,y] = h$, $[h,e_0] = e_0$, $[h,e_1] = -e_1$, $[x,e_1] = e_0$, $[y,e_0] = e_1$.\\
 $i = 1$, $c=3$, $p=1$, $F = L_5$, no CP's,\\
 $Y = k[f] = Sy$, where $f = e_0e_1h + e_1^2x - e_0^2y$\\
 $R^I = k(f)$, $M = k[e_0, e_1, f]$.
 \end{itemize}
 
 {\bf I.2. $\vet{\dim L = 6}$, with basis $\vet{h,x,y,e_0,e_1,e_2}$}
 \begin{itemize}
 \item[2.] $L_{6,1} = sl(2,k) \oplus W_2$ (quadratic)\\
 $[h,x] = 2x$, $[h,y] = -2y$, $[x,y] = h$, $[h,e_0] = 2e_0$, $[h,e_2] = -2e_2$,\\
 $[x,e_1] = 2e_0$, $[x,e_2] = e_1$, $[y,e_0] = e_1$, $[y,e_1] = 2e_2$\\
 $i = 2$, $c=4$, $p=1$ $F= L_{6,1}$, no CP's\\
 $Y = k[f_1,f_2] = Sy$, $f_1 = e_1^2 - 4e_0e_2$, $f_2 = e_1h + 2e_2x - 2e_0y$\\
$R^I = k(f_1, f_2)$, $M = k[e_0, e_1, e_2, f_2]$.
\item[3.] $L_{6,3} = sl(2,k) \oplus H$ (quasi quadratic) (see Example 54)\\
$[h,x] = 2x$, $[h,y] = -2y$, $[x,y] = h$, $[h,e_0] = e_0$, 
$[h,e_1] = -e_1$, $[x,e_1] = e_0$, $[y,e_0] = e_1$, $[e_0,e_1] = e_2$\\
 $i = 2$, $c=4$, $p=1$ $F= L_{6,3}$, no CP's\\
 $Y = k[e_2,f] = Sy$, $f = e_2(h^2 + 4xy) + 2(e_0e_1h + e_1^2x - e_0^2y)$, $R^I = k(e_2, f)$,\\
$M = k[e_1, e_2, e_0^2- 2e_2x, f]$.
\item[4.] $L_{6,4} = L_5  \oplus ke_2$ (Frobenius)\\
$[h,x] = 2x$, $[h,y] = -2y$, $[x,y] = h$, $[h,e_0] = e_0$, $[h,e_1] = -e_1$,\\
$[x,e_1] = e_0$, $[y,e_0] = e_1$, $[e_2, e_0] = e_0$, $[e_2, e_1] = e_1$.\\
$i = 0$, $c=3$, $p=e_0e_1h + e_1^2x - e_0^2y$, $F = 0$, no CP's,\\
$Y = k$, $Sy = k[p]$, $R^I = k$, $M = k[e_0, e_1, p]$.
\end{itemize}

{\bf I.3. $\vet{\dim L = 7}$, with basis $\vet{h, x, y, e_0, e_1, e_2, e_3}$}
\begin{itemize}
\item[5.] $L_{7,1} = sl(2,k) \oplus W_3$\\
$[h,x] = 2x$, $[h,y] = -2y$, $[x,y] = h$, $[h, e_0] = 3e_0$,\\
$[h, e_1] = e_1$, $[h,e_2] = -e_2$, $[h, e_3] = -3e_3$, $[x,e_1] = 3e_0$\\
$[x, e_2] = 2e_1$, $[x,e_3] = e_2$, $[y,e_0] = e_1$, $[y, e_1] = 2e_2$,\\
$[y, e_2] = 3e_3$.\\
$i = 1$, $c=4$, $p = 1$, $F = W_3 = CPI$\\
$Y = k[f] = k[W_3]^{SL(2)} = Sy$, $R^I = k(f)$,\\
$f = 4e_0e_2^3 - e_1^2e_2^2 - 18e_0e_1e_2e_3 + 27e_0^2e_3^2 + 4e_1^3e_3$\\
$M = k[e_0, e_1, e_2, e_3]$
\item[6.] $L_{7,2} = sl(2,k) \oplus W_1 \oplus W_1$\\
$[h,x] = 2x$, $[h,y] = -2y$, $[x,y] = h$, $[h, e_0] = e_0$,\\
$[h, e_1] = -e_1$, $[h,e_2] = e_2$, $[h, e_3] = -e_3$, $[x,e_1] = e_0$\\
$[x, e_3] = e_2$, $[y,e_0] = e_1$, $[y, e_2] = e_3$.\\
$i = 1$, $c=4$, $p = e_0e_3 - e_1e_2$, $F = W_1 \oplus W_1 = CPI$,\\
$Y = k[p] = Sy$, $R^I = k(p)$, $M = k[e_0, e_1, e_2, e_3]$.
\item[7.] $L_{7,7} = L_{6,1} \oplus ke_3$\\
$[h,x] = 2x$, $[h,y] = -2y$, $[x,y] = h$, $[h, e_0] = 2e_0$,\\
$[h, e_2] = -2e_2$, $[x,e_1] = 2e_0$, $[x,e_2] = e_1$, $[y, e_0] = e_1$,\\
$[y, e_1] = 2e_2$, $[e_3, e_0] = e_0$, $[e_3, e_1] = e_1$, $[e_3, e_2] = e_2$.\\
$i = 1$, $c=4$, $p = 1$, $F = L_{6,1} = (L_{7,7})_\Lambda$, no CP's.\\
$Y = k$, $Sy = k[f_1, f_2]$, $f_1 = e_1^2 - 4e_0e_2$, $f_2 = e_1h + 2e_2x - 2e_0y$,\\
$R^I = k(f_2^2/f_1)$, $M = k[e_0, e_1, e_2, f_2]$.
\item[8.] $L_{7,8} (\alpha \neq 0) = (L_5 \times ke_2) \oplus ke_3$\\
$[h,x] = 2x$, $[h,y] = -2y$, $[x,y] = h$, $[h, e_0] = e_0$,\\
$[h, e_1] = -e_1$, $[x,e_1] = e_0$, $[y, e_0] = e_1$, $[e_3,e_0] = e_0$\\
$[e_3, e_1] = e_1$, $[e_3,e_2] = \alpha e_2$\\
$i = 1$, $c=4$, $p = 1$, $F = L_5 \times ke_2 = (L_{7,8})_\Lambda$, no CP's\\
$Y = k$, $Sy = k[e_2,f]$, $f = e_0e_1h + e_1^2x - e_0^2y$,\\
$R^I = k(e_2^rf^s)$, $r, s \in \Bbb Z$ coprime such that $r\alpha + 2s = 0$,\\
$M = k[e_0, e_1, e_2, f]$.\\[0ex]
[If $\alpha \notin \Bbb Q$, then $R^I = k$]
\item[9.] $L_{7,9} = L_{6,3} \oplus ke_3$ (see Example 55)\\
$[h,x] = 2x$, $[h,y] = -2y$, $[x,y] = h$, $[h, e_0] = e_0$, $[h,e_1] = -e_1$,\\
$[x,e_1] = e_0$, $[y, e_0] = e_1$, $[e_0, e_1] = e_2$, $[e_3, e_0] = e_0$,\\
$[e_3, e_1] = e_1$, $[e_3, e_2] = 2e_2$.\\
$i=1$, $c=4$, $p=1$, $F = L_{6,3} = (L_{7,9})_\Lambda$, no CP's\\
$Y = k$, $Sy = k[e_2,f]$, $f = e_2(h^2 + 4xy) + 2(e_0e_1h + e_1^2x - e_0^2y)$\\
$R^I = k(f/e_2)$, $M = k[e_1, e_2, e_0^2-2e_2x, f]$
\end{itemize}
{\bf I.4. $\vet{\dim L = 8}$, with basis $\vet{h, x, y, e_0, e_1, e_2, e_3, e_4}$}
\begin{itemize}
\item[10.] $L_{8,1} = sl(2,k) \oplus W_4$\\
$[h,x] = 2x$, $[h,y] = -2y$, $[x,y] = h$, $[h, e_0] = 4e_0$, $[h,e_1] = 2e_1$,\\
$[h,e_3] = -2e_3$, $[h, e_4] = -4e_4$, $[x, e_1] = 4e_0$, $[x, e_2] = 3e_1$,\\
$[x, e_3] = 2e_2$, $[x, e_4] = e_3$, $[y,e_0] = e_1$, $[y,e_1] = 2e_2$, $[y,e_2] = 3e_3$, $[y,e_3] = 4e_4$.\\
$i=2$, $c=5$, $p=1$, $F = W_4 = CPI$,\\
$Y = k[f_1, f_2] = Sy$, $f_1 = e_2^2 - 3e_1e_3 + 12e_0e_4$,\\
$f_2 = 2e_2^3-9e_1e_2e_3 + 27e_0e_3^2 + 27e_1^2e_4 - 72 e_0e_2e_4$\\
$R^I = k(f_1, f_2)$, $M = k[e_0, e_1, e_2, e_3, e_4]$
\item[11.] $L_{8,2} = sl(2,k) \oplus W_2 \oplus W_1$\\
$[h,x] = 2x$, $[h,y] = -2y$, $[x,y] = h$, $[h, e_0] = 2e_0$, $[h,e_2] = -2e_2$,\\
$[h,e_3] = e_3$, $[h, e_4] = -e_4$, $[x, e_1] = 2e_0$, $[x, e_2] = e_1$,\\
$[x, e_4] = e_3$, $[y, e_0] = e_1$, $[y,e_1] = 2e_2$, $[y,e_3] = e_4$\\
$i=2$, $c=5$, $p=1$, $F = W_2 \oplus W_1 = CPI$,\\
$Y = k[f_1, f_2] = Sy$, $f_1 = e_1^2 - 4e_0e_2$, $f_2 = e_0e_4^2 - e_1e_3e_4 + e_2e_3^2$\\
$R^I = k(f_1, f_2)$, $M = k[e_0, e_1, e_2, e_3, e_4]$
\item[12.] $L_{8,13}$ (quasi quadratic)\\
$[h,x] = 2x$, $[h,y] = -2y$, $[x,y] = h$, $[h, e_0] = e_0$, $[h,e_1] = -e_1$,\\
$[h,e_3] = e_3$, $[h, e_4] = -e_4$, $[x, e_1] = e_0$, $[x, e_4] = e_3$,\\
$[y, e_0] = e_1$, $[y, e_3] = e_4$, $[e_0,e_1] = e_2$.
$i=2$, $c=5$, $p=1$, $F = L_{8,13}$, no CP's\\
$Y = k[e_2, f] = Sy$, $f = 2e_2 (e_3e_4h + e_4^2x - e_3^2y) - (e_0e_4 - e_1e_3)^2$\\
$R^I = k(e_2, f)$, $M = k[e_1, e_2, e_3, e_0^2 - 2e_2x, f]$
\item[13.] $L_{8,14} = (L_{6,3} \times ke_3) \oplus ke_4$\\
$[h,x] = 2x$, $[h,y] = -2y$, $[x,y] = h$, $[h, e_0] = e_0$, $[h,e_1] = -e_1$,\\
$[x,e_1] = e_0$, $[y, e_0] = e_1$, $[e_0, e_1] = e_2$, $[e_3, e_4] = e_2$.\\
$i=2$, $c=5$, $p=e_2$, $F = L_{6,3}$, no CP's\\
$Y = k[e_2, f] = Sy$, $f = e_2(h^2 + 4xy) + 2(e_0e_1h + e_1^2 x - e_0^2y)$,\\
$R^I = k(e_2, f)$, $M = k[e_1, e_2, e_3, e_0^2 - 2e_2x, f]$
\item[14.] $L_{8,15}$ (quasi quadratic)\\
$[h,x] = 2x$, $[h,y] = -2y$, $[x,y] = h$, $[h, e_0] = e_0$, $[h,e_1] = -e_1$,\\
$[h,e_3] = e_3$, $[h, e_4] = -e_4$, $[x, e_1] = e_0$, $[x, e_4] = e_3$,\\
$[y, e_0] = e_1$, $[y, e_3] = e_4$, $[e_0,e_1] = e_2$, $[e_3,e_4] = e_2$.\\
$i=2$, $c=5$, $p=1$, $F = L_{8,15}$, no CP's\\
$Y = k[e_2, f] = Sy$, $f = e_2^2(h^2 + 4xy) + 2e_2(e_0e_1 + e_3e_4) h + 2e_2 (e_1^2 + e_4^2)x -
2e_2 (e_0^2 + e_3^2)y - (e_0e_4 - e_1e_3)^2$,\\
$R^I = k(e_2, f)$, $M = k[e_1, e_2, e_4, e_0^2 + e_3^2 - 2e_2x, f]$.
\item[15.] $L_{8,16}$ (quasi quadratic)\\
$[h,x] = 2x$, $[h,y] = -2y$, $[x,y] = h$, $[h, e_1] = 3e_1$, $[h,e_2] = e_2$,\\
$[h,e_3] = -e_3$, $[h, e_4] = -3e_4$, $[x, e_2] = 3e_1$, $[x, e_3] = 2e_2$,\\
$[x, e_4] = e_3$, $[y, e_1] = e_2$, $[y,e_2] = 2e_3$, $[y,e_3] = 3e_4$, $[e_1, e_4] =e_0$, $[e_2, e_3] = -3e_0$.\\
$i=2$, $c=5$, $p=1$, $F = L_{8,16}$, no CP's\\
$Y = k[e_0, f] = Sy$, $f = 3e_0^2(h^2 + 4xy) + 2e_0(9e_1e_4 - e_2e_3) h + 4e_0 (3e_2e_4 - e_3^2)x +
4e_0 (e_2^2 - 3e_1e_3)y + 4e_1e_3^3 - e_2^2e_3^2 - 18e_1e_2e_3e_4 + 27 e_1^2 e_4^2 + 4e_2^3e_4$,\\
$R^I = k(e_0, f)$, $M = k[e_0, e_1, e_2, 3e_0x + e_2^2 - 3e_1e_3, f]$.
\item[16.] $L_{8,17} = L_{7,2} \oplus ke_4$ (quasi quadratic) (see Example 56)\\
$[h,x] = 2x$, $[h,y] = -2y$, $[x,y] = h$, $[h, e_0] = e_0$, $[h,e_1] = -e_1$,\\
$[h,e_2] = e_2$, $[h, e_3] = -e_3$, $[x, e_1] = e_0$, $[x, e_3] = e_2$,\\
$[y, e_0] = e_1$, $[y, e_2] = e_3$, $[e_2,e_4] = e_0$, $[e_3,e_4] = e_1$.\\
$i=2$, $c=5$, $p=1$, $F = L_{8,17}$, no CP's\\
$Y = k[f_1, f_2] = Sy$, $f_1 = e_1e_2 - e_0e_3$, $f_2 = e_0e_1h + e_1^2x - e_0^2y + (e_1e_2 - e_0e_3)e_4$\\
$R^I = k(f_1, f_2)$, $M = k[e_0, e_1, e_2, e_3, f_2]$.
\item[17.] $L_{8,18}$ (quasi quadratic)\\
$[h,x] = 2x$, $[h,y] = -2y$, $[x,y] = h$, $[h, e_0] = e_0$, $[h,e_1] = -e_1$,\\
$[h,e_3] = e_3$, $[h, e_4] = -e_4$, $[x, e_1] = e_0$, $[x, e_4] = e_3$,\\
$[y, e_0] = e_1$, $[y, e_3] = e_4$, $[e_2,e_3] = e_0$, $[e_2,e_4] = e_1$, $[e_3, e_4] = e_2$.\\
$i=2$, $c=5$, $p=1$, $F = L_{8,18}$, no CP's\\
$Y = k[f_1, f_2] = Sy$, $f_1 = 2(e_0e_4 - e_1e_3) + e_2^2$,\\
$f_2 = e_0e_1h + e_1^2x - e_0^2y + e_2(e_0e_4 - e_1e_3) + \ds\frac{1}{3}e_2^3$,\\
$R^I = k(f_1, f_2)$, $M = k[e_0, e_1, e_2, f_1, f_2]$.
\item[18.] $L_{8,19} = L_{7,1} \oplus ke_4$ (Frobenius)\\
$[h,x] = 2x$, $[h,y] = -2y$, $[x,y] = h$, $[h, e_0] = 3e_0$, $[h,e_1] = e_1$,\\
$[h,e_2] = -e_2$, $[h, e_3] = -3e_3$, $[x, e_1] = 3e_0$, $[x, e_2] = 2e_1$,\\
$[x, e_3] = e_2$, $[y, e_0] = e_1$, $[y,e_1] = 2e_2$, $[y,e_2] = 3e_3$, $[e_4, e_0] = e_0$, $[e_4, e_1] = e_1$, $[e_4,e_2] = e_2$, $[e_4, e_3] = e_3$.\\
$i=0$, $c=4$, $F=0$, $CPI = \langle e_0, e_1, e_2, e_3\rangle$, $(L_{8,19})_\Lambda = L_{7,1}$,\\
$p = 4e_0e_2^3 - e_1^2e_2^2 - 18e_0e_1e_2e_3 + 27e_0^2 e_3^2 + 4e_1^3e_3$,\\
$Y = k$, $Sy = k[p]$, $R^I = k$, $M = k[e_0, e_1, e_2, e_3]$.
\item[19.] $L_{8,20} (\alpha \neq -1) = L_{7,2} \oplus ke_4$ (Frobenius)\\
$[h,x] = 2x$, $[h,y] = -2y$, $[x,y] = h$, $[h, e_0] = e_0$, $[h,e_1] = -e_1$,\\
$[h,e_2] = e_2$, $[h, e_3] = -e_3$, $[x, e_1] = e_0$, $[x, e_3] = e_2$,\\
$[y, e_0] = e_1$, $[y, e_2] = e_3$, $[e_4,e_0] = e_0$, $[e_4,e_1] = e_1$, $[e_4, e_2] = \alpha e_2$, $[e_4, e_3] = \alpha e_3$.\\
$i=0$, $c=4$, $p= (e_0e_3 - e_1e_2)^2$, $F = 0$, $CPI = \langle e_0, e_1, e_2, e_3\rangle$, $(L_{8,20})_\Lambda = L_{7,2}$\\
$Y = k$, $Sy= k[e_0e_3 - e_1e_2]$, $R^I = k$, $M = k[e_0, e_1, e_2, e_3]$.
\item[20.] $L_{8,20} (\alpha = -1) = L_{7,2} \oplus ke_4$ (quasi quadratic)\\
$[h,x] = 2x$, $[h,y] = -2y$, $[x,y] = h$, $[h, e_0] = e_0$, $[h,e_1] = -e_1$,\\
$[h,e_2] = e_2$, $[h, e_3] = -e_3$, $[x, e_1] = e_0$, $[x, e_3] = e_2$,\\
$[y, e_0] = e_1$, $[y, e_2] = e_3$, $[e_4,e_0] = e_0$, $[e_4,e_1] = e_1$, $[e_4, e_2] = -e_2$, $[e_4, e_3] = -e_3$.\\
$i=2$, $c=5$, $p=1$, $F = L_{8,20}$, no CP's\\
$Y = k[f_1, f_2] = Sy$, $f_1 = e_1e_2 - e_0e_3$,\\
$f_2 = (e_1e_2 + e_0e_3)h + 2e_1e_3x - 2e_0e_2y + (e_1e_2 - e_0e_3)e_4$\\
$R^I = k(f_1, f_2)$, $M = k[e_0, e_1, e_2, e_3, f_2]$.
\item[21.] $L_{8,21} (\alpha \neq 0, \beta \neq 0) = (L_5 \times \langle e_2, e_3\rangle ) \oplus ke_4$\\
$[h,x] = 2x$, $[h,y] = -2y$, $[x,y] = h$, $[h, e_0] = e_0$, $[h,e_1] = -e_1$,\\
$[x,e_1] = e_0$, $[y, e_0] = e_1$, $[e_4, e_0] = e_0$, $[e_4, e_1] = e_1$,\\
$[e_4, e_2] = \alpha e_2$, $[e_4, e_3] = \beta e_3$.\\
$i=2$, $c=5$, $p=1$, $F = L_{5}\times \langle e_2, e_3\rangle = (L_{8,21})_\Lambda$, no CP's\\
$Y = k$, $Sy = k[e_2, e_3,f]$, $f = e_0e_1h + e_1^2x - e_0^2y$, \\
$M = k[e_0, e_1, e_2, e_3, f]$, $ R^I = k(e_2^{r_0} e_3^{s_o} f^{t_0}, e_2^{r_1} e_3^{s_1} f^{t_1})$, where $(r_0, s_0, t_0)$, $(r_1, s_1, t_1)$ is a basis of the free
$\Bbb Z$-module $\{(r, s, t) \in \Bbb Z^3 \mid \alpha r + \beta s + 2t = 0\}$.\\[0ex]
[If $\alpha \notin \Bbb Q$ then $R^I = k(e_3^s f^t)$ where $s,t$ are coprime integers such that $\beta s + 2t = 0$.
Similarly, if $\beta \notin \Bbb Q$.  If $\alpha$ and $\beta \notin \Bbb Q$ and $\alpha / \beta \notin \Bbb Q$ then $R^I = k$]
\item[22.] $L_{8,22} (\alpha \neq 0) = (L_{6,1} \times ke_3) \oplus ke_4$\\
$[h,x] = 2x$, $[h,y] = -2y$, $[x,y] = h$, $[h, e_0] = 2e_0$, $[h,e_2] = -2e_2$,\\
$[x,e_1] = 2e_0$, $[x, e_2] = e_1$, $[y, e_0] = e_1$, $[y, e_1] = 2e_2$,\\
$[e_4, e_0] = e_0$, $[e_4, e_1] = e_1$, $[e_4, e_2] = e_2$, $[e_4, e_3] = \alpha e_3$.\\
$i=2$, $c=5$, $p=1$, $F = L_{6,1}\times ke_3 = (L_{8,22})_\Lambda$, no CP's\\
$Y = k$, $Sy = k[e_3,f_1, f_2]$, $f_1 = e_1^2 - 4e_0e_2$, \\
$f_2 = e_1h + 2e_2x - 2e_0y$, $M = k[e_0, e_1, e_2, e_3, f_2]$,\\
$R^I = k(f_2^2/f_1, f_2^s / e_3^t)$ where $s$, $t$ are coprime integers such that $\alpha = \ds\frac{s}{t}$.\\[0ex]
[If $\alpha \notin \Bbb Q$ then $R^I = k(f_2^2 / f_1)$]
\item[23.] $L_{8,23} (\alpha \neq 0) = (L_{6,3} \times ke_3) \oplus ke_4$\\
$[h,x] = 2x$, $[h,y] = -2y$, $[x,y] = h$, $[h, e_0] = e_0$, $[h,e_1] = -e_1$,\\
$[x,e_1] = e_0$, $[y, e_0] = e_1$, $[e_0, e_1] = e_2$, $[e_4, e_0] = e_0$,\\
$[e_4, e_1] = e_1$, $[e_4, e_2] = 2e_2$, $[e_4, e_3] = \alpha e_3$.\\
$i=2$, $c=5$, $p=1$, $F = L_{6,3}\times ke_3 = (L_{8,23})_\Lambda$, no CP's\\
$Y = k$, $Sy = k[e_2,e_3, f]$, $f = e_2(h^2 + 4xy) + 2(e_0e_1h + e_1^2x - e_0^2y)$,\\
$M = k[e_1, e_2, e_3, e_0^2 - 2e_2x, f]$,\\
$R^I = k(f/e_2, e_2^se_3^t)$ where $s$, $t$ are coprime integers such that $2s + \alpha t = 0$.\\[0ex]
[If $\alpha \notin \Bbb Q$ then $R^I = k(f / e_2)$]
\item[24.] $L_{8,24} = sl(3,k)$ (simple and hence quadratic)\\
Basis: (see e.g. [D4])\\
$h_\alpha = E_{11} - E_{22}$, $h_\gamma = E_{22} - E_{33}$, $x_\alpha = E_{12}$, $x_\beta = E_{13}$, $x_\gamma = E_{23}$,
$x_{-\alpha} = E_{21}$, $x_{-\beta} = E_{31}$, $x_{-\gamma} = E_{32}$, where the $E_{ij}$ are the standard $3 \times 3$ matrices and $\alpha, \gamma, \beta = \alpha + \gamma$ are the
positive roots w.r.t. the Cartan subalgebra $H = \langle E_{11} - E_{22}, E_{22} - E_{33}\rangle$.\\
$[h_\alpha, x_\alpha] = 2x_\alpha$, $[h_\alpha, x_\beta] = x_\beta$, $[h_\alpha, x_\gamma] = -x_\gamma$, $[h_\alpha, x_{-\alpha}] = -2x_{-\alpha}$,
$[h_\alpha, x_{-\beta}] = -x_{-\beta}$, $[h_\alpha, x_{-\gamma}] = x_{-\gamma}$, $[h_\gamma, x_\alpha] = -x_\alpha$, $[h_\gamma, x_\beta] = x_\beta$,
$[h_\gamma, x_\gamma] = 2x_\gamma$, $[h_\gamma, x_{-\alpha}] = x_{-\alpha}$, $[h_\gamma, x_{-\beta}] = -x_{-\beta}$, $[h_\gamma, x_{-\gamma}] = -2x_{-\gamma}$,
$[x_\alpha, x_\gamma] = x_\beta$, $[x_\alpha, x_{-\alpha}] = h_\alpha$, $[x_\alpha, x_{-\beta}] = -x_{-\gamma}$, $[x_\beta, x_{-\alpha}] = -x_\gamma$,
$[x_\beta, x_{-\beta}] = h_\alpha + h_\gamma$, $[x_\beta, x_{-\gamma}] = x_{\alpha}$, $[x_\gamma, x_{-\beta}] = x_{-\alpha}$,
$[x_\gamma, x_{-\gamma}] = h_\gamma$, $[x_{-\alpha}, x_{-\gamma}] = -x_{-\beta}$\\
$i=2$, $c=5$, $p=1$, $F = L_{8,24}$, no CP's,\\
$Y = k[f_1,f_2] = Sy$,\\
$f_1 = h_\alpha^2 + h_\alpha h_\gamma + h_\gamma^2 + 3(x_\alpha x_{-\alpha} + x_\beta x_{-\beta} + x_\gamma x_{-\gamma})$\\
$f_2 = (h_\alpha + 2h_\gamma) (h_\alpha - h_\gamma) (2h_\alpha + h_\gamma) + 9(h_\alpha + 2h_\gamma) x_\alpha x_{-\alpha} + \\
9 (h_\alpha - h_\gamma) x_\beta x_{-\beta} - 9(2h_\alpha + h_\gamma) x_\gamma x_{-\gamma} + 27 (x_\alpha x_\gamma x_{-\beta} + x_\beta x_{-\alpha} x_{-\gamma})$\\
$R^I = k(f_1, f_2)$, $M = k[f_1, f_2, f_3, f_4, f_5]$, $f_3 = x_\beta$,\\
$f_4 = x_\alpha + x_\gamma$, $f_5 = (h_\alpha + 2h_\gamma) x_\alpha - (2h_\alpha + h_\gamma) x_\gamma + 3x_\beta (x_{-\alpha} + x_{-\gamma})$.
\end{itemize}
The subalgebra $M$ was constructed by means of the argument shift method starting out from the generating invariants $f_1, f_2$ of $Y(L_{8,24})$.  By [PY1] $M$ is a polynomial,
strongly complete Poisson commutative subalgebra of $S(L_{8,24})$.  See also [Ta].\\
\ \\
{\bf II. $\vet{L}$ is not algebraic}\\
This includes the families of part I with non rational parameters.\\
In the remaining cases $L$ will be 8-dimensional with basis $h, x, y, e_0, e_1, e_2, e_3, e_4$.
\begin{itemize}
\item[25.] $L_{8,25} = (L_5 \times \langle e_2, e_3\rangle )\oplus ke_4$ (see Example 58)\\
$[h,x] = 2x$, $[h,y] = -2y$, $[x,y] = h$, $[h, e_0] = e_0$, $[h,e_1] = -e_1$,\\
$[x,e_1] = e_0$, $[y, e_0] = e_1$, $[e_4, e_0] = e_0$, $[e_4, e_1] = e_1$, $[e_4, e_2] = -e_3$\\
$i=2$, $c=5$, $p=1$, $F = L_5 \times \langle e_2, e_3\rangle = (L_{8,25})_\Lambda$, no CP's\\
$Y = k[e_3]$, $Sy = k[e_3, f]$, $f = e_0e_1 h + e_1^2x - e_0^2y$,\\
$R^I = k(e_3)$, $M = k[e_0, e_1, e_2, e_3, f]$.
\item[26.] $L_{8,26}(\alpha \neq 0) = (L_5 \times \langle e_2, e_3\rangle )\oplus ke_4$\\
$[h,x] = 2x$, $[h,y] = -2y$, $[x,y] = h$, $[h, e_0] = e_0$, $[h,e_1] = -e_1$,\\
$[x,e_1] = e_0$, $[y, e_0] = e_1$, $[e_4, e_0] = \alpha e_0$, $[e_4, e_1] = \alpha e_1$, $[e_4, e_2] = e_2$,\\
$[e_4, e_3] = e_3-e_2$.\\
$i=2$, $c=5$, $p=1$, $F = L_5 \times \langle e_2, e_3\rangle = (L_{8,26})_\Lambda$, no CP's\\
$Y = k$, $Sy(L) = k[e_2, f]$, $f = e_0e_1 h + e_1^2x - e_0^2y$,\\
$M = k[e_0, e_1, e_2, e_3, f]$
\begin{itemize}
\item[(i)] If $\alpha \in \Bbb Q$ then $R^I = k(f^t/e_2^s)$ where $s, t$ are coprime integers such that $\ds\frac{s}{t} = 2\alpha$.
\item[(ii)] If $\alpha \notin \Bbb Q$ then $R^I = k$.
\end{itemize}
\item[27.] $L_{8,26} (\alpha = 0) = (L_5 \times \langle e_2, e_3\rangle )\oplus ke_4$\\
$[h,x] = 2x$, $[h,y] = -2y$, $[x,y] = h$, $[h, e_0] = e_0$, $[h,e_1] = -e_1$,\\
$[x,e_1] = e_0$, $[y, e_0] = e_1$, $[e_4, e_2] = e_2$, $[e_4, e_3] = e_3-e_2$.\\
$i=2$, $c=5$, $p=1$, $F = L_5 \times \langle e_2, e_3\rangle = (L_{8,26})_\Lambda$, no CP's\\
$Y = k[f]$, $f = e_0e_1 h + e_1^2x - e_0^2y$, $ Sy = k[e_2, f]$,\\
$R^I = k(f)$, $M = k[e_0, e_1, e_2, e_3, f]$.
\item[28.] $L_{8,27} = (L_{6,3} \times ke_3) \oplus ke_4$\\
$[h,x] = 2x$, $[h,y] = -2y$, $[x,y] = h$, $[h, e_0] = e_0$, $[h,e_1] = -e_1$,\\
$[x,e_1] = e_0$, $[y, e_0] = e_1$, $[e_0, e_1] = e_2$, $[e_4, e_0] = e_0$, $[e_4, e_1] = e_1$,\\
$[e_4, e_2] = 2e_2$, $[e_4, e_3] = 2e_3-e_2$.\\
$i=2$, $c=5$, $p=1$, $F = L_{6,3} \times ke_3 = (L_{8,27})_\Lambda$, no CP's\\
$Y = k$, $Sy = k[e_2,f]$, $f = e_2(h^2 + 4xy) + 2(e_0e_1 h + e_1^2x - e_0^2y)$,\\
$R^I = k(f/e_2)$, $M = [e_1, e_2, e_3, e_0^2 - 2e_2x, f]$.
\item[29.] $L_{8,28} = L_{7,2} \oplus ke_4$ (Frobenius)\\
$[h,x] = 2x$, $[h,y] = -2y$, $[x,y] = h$, $[h, e_0] = e_0$, $[h,e_1] = -e_1$,\\
$[h,e_2] = e_2$, $[h, e_3] = -e_3$, $[x, e_1] = e_0$, $[x, e_3] = e_2$, $[y, e_0] = e_1$, $[y, e_2] = e_3$, $[e_4, e_0] = e_0$,
$[e_4, e_1] = e_1$, $[e_4, e_2] = e_2-e_0$, $[e_4, e_3] = e_3 - e_1$.\\
$i=0$, $c=4$, $p=(e_0e_3 - e_1e_2)^2$, $F = 0$, $\langle e_0, e_1, e_2, e_3\rangle = CPI$\\
$Y = k$, $Sy = k[e_0e_3 - e_1e_2]$, $R^I = k$, $M = k[e_0, e_1, e_2, e_3]$.
\end{itemize}

{\bf 5.2. Counterexample in dimension 9}\\
\ \\
{\bf Example 59.} Take the semi-direct product $L = sl(2,k) \oplus W_5$ in which $W_5$ is an abelian ideal.  As $[L,L] = L$, $L$ is algebraic without proper semi-invariants.  
Since $\dim sl(2,k) < \dim W_5$ we know that the stabilizer $sl(2,k)(f) = 0$ for some $f \in W_5^\ast$ by [AVE].  This implies that $i(L) = \dim W_5 - \dim sl(2,k) = 
3$, $Y(L) = S(W_5)^{sl(2,k)}$
and also that $W_5$ is a CPI of $L$ by [O4, Proposition 17].  One verifies that $\codim\ L_{\sing}^\ast = 4$.  By Proposition 15 (or by Theorem 16) we may conclude that $L$ is not coregular.
On the other hand, $L$ satisfies the Gelfand-Kirillov conjecture by [O6, Proposition 4.3].\\
\ \\
Now suppose $k = \Bbb C$.  Then $W_5$ may be considered as the vector space of binary forms of degree 5 with complex coefficients (the quintics) on which $SL(2,\Bbb C)$ acts.
The algebra of invariants $\Bbb C[W_5]^{SL(2,\Bbb C)}$ which is isomorphic to $S(W_5)^{sl(2,\Bbb C)}$, has been studied already in the 19th century by Sylvester, among others.  
At first 3 algebraically independent
invariants $I_4$, $I_8$, $I_{12}$ were found of degrees 4, 8, 12.  In 1854 Hermite discovered an invariant $I_{18}$ of degree 18 and he showed that
$$\Bbb C[W_5]^{SL(2,\Bbb C)} = \Bbb C[I_4, I_8, I_{12}, I_{18}]$$
with the following relation:
$$16I_{18}^2 = I_4I_8^4 + 8I_8^3I_{12} - 2I_4^2I_8^2I_{12} - 72I_4I_8I_{12}^2 - 432I_{12}^3 + I_4^3 I_{12}^2$$
In particular, $\Bbb C[W_5]^{SL(2,\Bbb C)}$ is not polynomial.\\
The explicit forms of $I_4$, $I_8$, $I_{12}$, $I_{18}$ were given in papers by Cayley.  $I_{18}$ has 848 monomials with very large coefficients ! See [D5, p.41].\\
\ \\
{\bf Example 60.}\\
We conclude this section by considering the semi-direct product\\
$L = sl(2,k)\oplus W_2 \oplus W_2$ with standard basis $h,x,y$; $e_0,e_1,e_2$; $e_3,e_4,e_5$.\\
We know that this is a counterexample to the Gelfand-Kirillov conjecture [AOV1].  However, besides this, its behaviour is rather tame.\\
Indeed, $[L,L] = L$ so $L$ is algebraic and it has no proper semi-invariants.  Also, $F(L) = W_2 \oplus W_2$ is a CPI of $L$, $i(L) = 3$, $c(L) = 6$ and $p_L = 1$.\\
Furthermore, $L$ is coregular as\\
$Y(L) = k[f_1,f_2,f_3] = Sy(L)$, where $f_1 = e_1^2 - 4e_0e_2$, $f_2 = e_4^2 - 4e_3e_5$,\\
$f_3 = e_1e_4 - 2e_2e_3 - 2e_0e_5$ by Theorem 29.  
Consequently, $R(L)^L = k(f_1, f_2, f_3)$.\\
\ \\
{\bf 6. Dixmier's fourth problem}\\
\ \\
Let $L$ be a finite dimensional Lie algebra over an algebraically closed field $k$ of characteristic zero.  Then we know that the field $Z(D(L))$ is isomorphic with $R(L)^L$ and hence is an extension of finite type of $k$. [RV, p.401], [D6, 10.5.6].\\
In his book Enveloping Algebras Dixmier raised the following problem [D6, p.354] and proved it for $L$ solvable (In fact he even showed it for $L$ completely solvable over an arbitrary field $k$ of characteristic zero
[D6, Proposition 4.4.8].  It also holds for solvable $L$ over $k = \IR$ [Be]).\\
\ \\
{\bf Problem 61.}\\
Is $Z(D(L))$ rational over $k$ ? (i.e. is it a purely transcendental extension of $k$ ?).\\
To our knowledge this problem is still open.  Notice that the Gelfand-Kirillov conjecture is a much stronger condition.\\
Obviously, Dixmier's question has a positive answer if $L$ is coregular without proper semi-invariants, since then $Z(D(L))$ is precisely the quotient field of $Z(U(L))$, the latter being polynomial.  
This is especially the case for $L$ semi-simple and also for the canonical truncation $\mathfrak{g}_\Lambda$ of a Frobenius Lie algebra $\mathfrak{g}$ (since $Z(U(\mathfrak{g}_\Lambda)) = Sz (U(\mathfrak{g}))$,
which is polynomial [DNO]). Dixmier's statement is also true for a Lie algebra $L$ for which $j(L) = \dim Z(L)$ (Theorem 37).\\
\ \\
{\bf Proposition 62.} Consider the semi-direct product $L = sl(2,k) \oplus W_n$, where $W_n$ is the $(n+1)$-dimensional irreducible $sl(2,k)$-module.
Then $Z(D(L))$ is rational over k.\\
\ \\
On the other hand, $L$ does not satisfy the Gelfand-Kirillov conjecture if $n$ is even and $n \geq 6$. [O6, Proposition 4.3].\\
\ \\
{\bf Proof.} We may assume that $n \geq 5$, since we verified it for $n=1, 2, 3, 4$ (see $L_5$, $L_{6,1}$, $L_{7,1}$, $L_{8,1}$ of 5.1).  Hence, $\dim sl(2,k) = 3 < n+1 = \dim W_n$.
We now follow the same argument as in Example 59.  By [AVE]
$$sl(2,k)(f) = 0\ \mbox{for some}\ f \in W_n^\ast$$
which by [O4, Proposition 17] implies that $i(L) = \dim W_n - \dim sl(2,k) = n+1 - 3 = n-2$ and
$$Z(D(L)) = R(W_n)^{sl(2,k)} = R(W_n)^{SL(2,k)}$$
where the latter is rational over $k$ [BK]. \hfill $\square$\\
\ \\
Next, Dixmier's result above combined with Theorem 53 yields:\\
(since the rationality of $R(L)^L$ is preserved under taking direct products)\\
\ \\
{\bf Proposition 63.} Assume $L$ is a Lie algebra over $k$ of dimension at most 8.  Then $Z(D(L))$ is rational over $k$.\\
\ \\
{\bf Lemma 64.} Let $L$ be an algebraic Lie algebra.  Then there exists a torus $T \subset L$ such that $L = L_\Lambda \oplus T$ and a basis $t_1, \ldots, t_r$ of $T$ such that 
$\ad\ t_1, \ldots, \ad\ t_r$ have rational eigenvalues and such that $\lambda(t_i) \in \Bbb Q$ for all $\lambda \in \Lambda(L)$, $i = 1, \ldots, r$.\\
\ \\
{\bf Proof.} As $L$ is algebraic, so is $\ad\ L$. By [C, p.324] $L$ admits the following decomposition
$$L = S \oplus N \oplus A$$
where $S$ is a semi-simple Lie subalgebra of $L$, $N$ is the nilradical of $L$ and $A$ is a torus of $L$ (i.e. an abelian Lie subalgebra of $L$ such that $\ad_LA$ consists of semi-simple elements).  
Also, $[A,S] = 0$,
$R = N \oplus A$ is the (solvable) radical of $L$ and $\ad_LA$ is algebraic.  It follows that also $R$ is algebraic [C, p.309].\\
By [Mc, Theorem 3.3] there is a basis $a_1, \ldots, a_s$ of $A$ such that $\ad_Na_1,\ldots,  \ad_Na_s$ have rational eigenvalues.  The same holds for $ad_La_1, \ldots, ad_La_s$ since $A$ is 
abelian and $[A,S] = 0$.  Because $S = [S,S]$ and each $x \in N$ acts locally nilpotent on $U(L)$, we see that $S \oplus N \subset L_\Lambda$.  Hence, $L_\Lambda + A = L$.  We may
assume that $L_\Lambda \neq L$ (otherwise take $T=0$).  Let $t_1$ be the first one among $a_1,\ldots, a_s$ such that $t_1 \notin L_\Lambda$. Next, consider $L_\Lambda \oplus kt_1 \subset L$ and so on.  After
a number of steps we obtain
$$L_\Lambda \oplus \langle t_1,\ldots, t_r\rangle = L$$
where $t_1,\ldots, t_r$ are linearly independent over $k$ and $\ad_Lt_1,\ldots, \ad_Lt_r$ have rational eigenvalues.  So it suffices to put $T = \langle t_1,\ldots, t_r\rangle$.
Let $x_1,\ldots, x_n$ be a basis of $L$ such that for all $i = 1, \ldots, r; j=1,\ldots, n: \ad\ t_i(x_j) = q_{ij} x_j$ for some $q_{ij} \in \Bbb Q$.  Next, let $u \in U(L)$ be a nonzero semi-invariant with weight $\lambda$, which can be written as 
$$u = \sum\limits_m \alpha_mx_1^{m_1} \ldots x_n^{m_n}\ \mbox{for some}\ \alpha_m \in k$$
and $m = (m_1,\ldots, m_n)$.  Now observe that
$$\sum\limits_m \lambda (t_i) \alpha_m x_1^{m_1} \ldots x_n^{m_n} = \lambda (t_i) u = \ad\ t_i(u)
= \sum\limits_m \left(\sum\limits_{j=1}^n m_j q_{ij}\right) \alpha_m x_1^{m_1} \ldots x_n^{m_n}$$
Now, select $m = (m_1,\ldots, m_n)$ such that $\alpha_m \neq 0$.  Then we may conclude that $\lambda(t_i) = \sum\limits_{j=1}^n m_j q_{ij} \in \Bbb Q$, 
$i = 1,\ldots, r$. \hfill $\square$\\
\ \\
{\bf Lemma 65.}
\begin{itemize}
\item[(1)] $Z(D(L)) \subset Z(D(L_\Lambda))$
\item[(2)] Assume $L$ is almost algebraic.  Then the field $Z(D(L_\Lambda))$ is generated by the semi-invariants of $U(L)$.
\item[(3)] If $L$ is algebraic then
$$\trdeg_kZ(D(L_\Lambda)) - \trdeg_k(Z(D(L)) = \dim L - \dim L_\Lambda$$
\end{itemize}

{\bf Proof.}  We recall that $Sz(U(L)) \subset Z(U(L_\Lambda))$ and equality occurs if $L$ is almost algebraic [DNO, Theorem 1.19].
\begin{itemize}
\item[(1)] Because any  nonzero $z \in Z(D(L))$ can be written as $z = uv^{-1}$, where $u$, $v$ are nonzero semi-invariants of $U(L)$ with the same weight [RV, Th\'eor\`eme 4.4],
[DNO, Corollary 1.10], we deduce at once that 
$$Z(D(L)) \subset Q(Sz(U(L))) \subset Q(Z(U(L_\Lambda))) = Z(D(L_\Lambda))$$
the latter since $L_\Lambda$ has no proper semi-invariants by 2 of Theorem 4.
\item[(2)] By definition $Sz(U(L))$ is generated by the semi-invariants of $U(L)$.
Hence, the same holds for its quotient field $Q(Sz(U(L))) = Q(Z(U(L_\Lambda))) = Z(D(L_\Lambda))$.
\item[(3)] If $L$ is algebraic then so is $L_\Lambda$ [DNO, Proposition 1.14].  By 4 of Theorem 4 we know that $c(L_\Lambda) = c(L)$.  Therefore
$$\dim L_\Lambda + i(L_\Lambda) = 2c(L_\Lambda) = 2c(L) = \dim L + i(L)$$
and hence $i(L_\Lambda) - i(L) = \dim L - \dim L_\Lambda$.\\
On the other hand,
$$i(L) = \trdeg_kZ(D(L))\ \ \ \mbox{and}\ \ \ i(L_\Lambda) = \trdeg_kZ(D(L_\Lambda))$$
by Theorem 1.\hfill $\square$
\end{itemize}
The following is the main result of this section.  It proved to be a useful tool in obtaining the explicit description of $R(L)^L$ in the list of 5.1.\\
\ \\
{\bf Theorem 66.} Let $L$ be an algebraic Lie algebra for which the field $Z(D(L_\Lambda))$ is freely generated by semi-invariants $u_1,\ldots, u_s$ of $U(L)$.  Then $Z(D(L))$ (and also $R(L)^L$) is rational over $k$.\\
\ \\
{\bf Proof.} By Lemma 64 there is a torus $T \subset L$ such that $L = L_\Lambda \oplus T$ and a basis $t_1,\ldots, t_r$ of $T$ such that for all $i = 1,\ldots, r$; 
$j = 1,\ldots, s$:
$$\ad\ t_i (u_j) = a_{ij} u_j \ \mbox{for some}\ a_{ij} \in \Bbb Q \hspace*{2cm} (\ast)$$
We may assume that $a_{ij} \in \Bbb Z$, since we can replace $t_i$ by a suitable integer multiple of itself.
Put $A = (a_{ij}) \in \Bbb Z^{r\times s}$ and $d_i = \ad\ t_i$.  For any $m = (m_1,\ldots, m_s) \in \Bbb Z^s$ we set
$$u^m = u_1^{m_1} \ldots u_s^{m_s}$$
Clearly, for all $i = 1, \ldots r$
$$d_i(u^m) = \ad\ t_i(u^m) = \left(\sum\limits_{j=1}^s a_{ij} m_j\right) u^m = (Am)_i u^m$$
Next, we consider
\begin{eqnarray*}
K &=& \left\{ m = (m_1,\ldots, m_s) \in \Bbb Z^s \mid d_i(u^m) = 0, i = 1,\ldots, r\right\}\\
&=& \left\{m \in \Bbb Z^s \mid Am = 0\right\}
\end{eqnarray*}
$K$ is a free $\Bbb Z$-module, being a submodule of the free $\Bbb Z$-module $\Bbb Z^s$.\\
Let $k_1, \ldots, k_q$ be a basis of $K$, then
$$z_1 = u^{k_1},\ldots, z_q=u^{k_q} \in Z(D(L))$$
(since they are annihilated by both $\ad\ L_\Lambda$ and $\ad\ T$)\\
\ \\
{\bf Claim:}\\
(1) $z_1,\ldots, z_q$ are algebraically independent over $k$.\\
(2) $Z(D(L)) = k(z_1,\ldots, z_q)$
\begin{itemize}
\item[(1)] By assumption $u_1,\ldots, u_s$ are algebraically independent over $k$, which is equivalent with the fact that \\
$u^m = u_1^{m_1} \ldots u_s^{m_s}$, $m = (m_1,\ldots, m_s) \in \Bbb N^s$ (and even $m \in \Bbb Z^s$) are linearly independent over $k$ ($\ast\ast$).  
Similarly, we have to show that $z_1^{m_1} \ldots z_q^{m_q}$, $m = (m_1,\ldots, m_q) \in \Bbb N^q$ are linearly independent over $k$.  Clearly,
$$z_1^{m_1} \ldots z_q^{m_q} = (u^{k_1})^{m_1} \ldots (u^{k_q})^{m_q} = u^{\sum m_ik_i}$$
and these are indeed linearly independent over $k$ by ($\ast\ast$), since $k_1,\ldots, k_q$ are linearly independent over $\Bbb Z$.
\item[(2)]We already know that $k(z_1,\ldots, z_q) \subset Z(D(L))$.  On the other hand, take $0 \neq z \in Z(D(L))$.  Then $z \in Z(D(L_\Lambda)) = k(u_1,\ldots, u_s)$ by (1) of the
previous lemma.  So, $z = vw^{-1}$, where $v$, $w$ are nonzero coprime elements of the polynomial algebra $k[u_1,\ldots, u_s]$ in the variables $u_1, \ldots, u_s$.  
We consider the degree of $v$ and $w$ with respect to these variables. Obviously $vw = wv$ and thus $z = vw^{-1} = w^{-1}v$.  Because $z \in Z(D(L))$ it is annihilated by each derivation $d_i = \ad\ t_i$, $i = 1,\ldots, r$.
 Therefore,
 $$zd_i(w) = d_i(z)w + zd_i(w) = d_i(zw) = d_i(v)$$
 and so $vd_i(w) = wd_i(v)$.  This implies that $v$ divides $d_i(v)$ as $v$ and $w$ are coprime. But $\deg(d_i(v)) \leq \deg (v)$ (by ($\ast$)) and thus $d_i(v) = \lambda_iv$ for a
 suitable $\lambda_i \in k$.  It follows that $d_i(w) = \lambda_iw$.  We can find nonzero $a_m \in k$ and $M \subset \Bbb N^s$ such that $v = \sum\limits_{m \in M} a_m u^m$.  Similarly,
$w =\sum\limits_{n \in N} b_n u^n$.  From
$$\sum\limits_{m \in M} a_m(Am)_i u^m = \sum\limits_{m \in M} a_m d_i(u^m) = d_i(v) = \lambda_iv = \sum\limits_{m \in M} a_m \lambda_i u^m$$ 
we obtain $(Am)_i = \lambda_i$, $m \in M$, $i=1,\ldots, r$.\\
Similarly, $(An)_i = \lambda_i$, $n \in N$, $i=1,\ldots, r$. Hence,\\
$(A(n-m))_i = (An)_i - (Am)_i = 0$.\\
Consequently, $A(n-m) = 0$ and thus $n-m \in K$ for all $n \in N$, $m \in M$.\\
We can find $\alpha_i \in \Bbb Z$ such that $n-m = \sum\limits_{i=1}^q \alpha_i k_i$.  Therefore,
$$u^{n-m} = u^{\sum \alpha_i k_i} =(u^{k_1})^{\alpha_1} \ldots (u^{k_q})^{\alpha_q}
= z_1^{\alpha_1} \ldots z_q^{\alpha_q} \in k (z_1,\ldots, z_q)$$
\end{itemize}
Finally,
\begin{eqnarray*}
z = vw^{-1} &=& \left(\sum\limits_m a_mu^m\right)\left(\sum\limits_n b_n u^n\right)^{-1} = \sum\limits_m a_m \left(u^{-m} \sum\limits_n b_n u^n\right)^{-1}\\
&=& \sum\limits_m a_m \left(\sum\limits_n b_n u^{n-m}\right)^{-1} \in k(z_1,\ldots, z_q)
\end{eqnarray*}
This establishes the claim, i.e. $Z(D(L))$ is rational over $k$.\hfill $\square$\\
\ \\
The following result by Panyushev [Pa1] can now be derived from Joseph's work on biparabolics.\\
\ \\
{\bf Corollary 67.} Let $L$ be a biparabolic (seaweed) subalgebra of a simple Lie algebra of type $A$ or $C$.  Then $R(L)^L$ (and hence also $Z(D(L)))$ is rational over $k$.\\
\ \\
{\bf Proof.} $L$ is algebraic [F, 6.4].  By [J3, J4] $Y(L_\Lambda)$ is freely generated by some semi-invariants of $S(L)$, say $v_1, \ldots, v_s$.  
By the Duflo isomorphism the same holds for $Z(U(L_\Lambda))$ and also for its quotient field, which is $Z(D(L_\Lambda))$.  By Theorem 66 $Z(D(L))$ is rational over $k$.\\
\ \\
{\bf Remark 68.} Most (but not all [Y]) (bi)parabolic subalgebras of semi-simple Lie algebras have canonical truncations whose Poisson centers are freely generated by semi-invariants
[F, FJ, FJ2, J3, J4].  Hence, we may then draw the same conclusion as above.\\
Let $L$ be a finite dimensional algebraic Lie algebra over $k$ and $G$ its algebraic adjoint group.  $G$ and $L$ act on $L^\ast$ via the coadjoint action.  We identify $R(L)$ with the field of rational
functions on $L^\ast$.\\
\ \\
{\bf Definition 69.} (See e.g. [TY1]) An affine slice of $L$ is an affine subspace $V$ of $L^\ast$ such that there exists an open subset $U$ of $V$ verifying the following conditions:
\begin{itemize}
\item[(1)] The set $G.U$ is dense in $L^\ast$
\item[(2)] $T_f(G.f) \cap T_f(U) = \{0\}$ for all $f \in U$
\item[(3)] $G.f \cap U = \{f\}$ for all $f \in U$
\end{itemize}
Such an affine slice exists for the coadjoint action for certain truncated [J6, J7] and non truncated [TY1] biparabolic subalgebras of a semi-simple Lie algebra.\\
\ \\
{\bf Theorem 70.} (Tauvel, Yu [TY1, Theorem 3.3.1])\\
Let $L$ be algebraic with algebraic adjoint group $G$.  Suppose there exists an affine slice for the coadjoint action of $L$.  Then $R(L)^G$ is rational aver $k$.  Hence the same holds for $Z(D(L))$ (since
the latter is isomorphic to $R(L)^G$ by [RV, 4.5]).\\
\ \\

{\bf 7. Appendix}\\
\ \\
{\bf Example 71.} Let $L$ be the nonalgebraic Lie algebra over $\Bbb C$ with basis $x, y, z, t$ and nonzero brackets.
$$[x,y] = y, \ \ [x,z] = \alpha z, \ \ [y,z] = t, \ \ [x,t] = (1+\alpha) t$$
with $\alpha$ irrational.\\
This example was introduced in [GK, p.522] in order to demonstrate the existence of nonalgebraic Lie algebras satisfying the Gelfand-Kirillov conjecture.  However the proof is incorrect, as the given Weyl generators of
$D(L)$, namely
$$p_1 = yt^{-1}, \ \ q_1 = z, \ \ p_2 = (1+\alpha)^{-1} t, \ \ q_2 = yzt^{-2} x t^{-1}$$
do not satisfy the necessary requirements, for instance\\
$[q_1, q_2] = -((x + 1)t + \alpha yz)zt^{-3} \neq 0$.  Probably a term is missing in $q_2$.\\
We now present a very short proof.  First we observe that $L$ can be considered as the semi-direct product $\mathfrak{g} \oplus W$ of the Lie algebra $\mathfrak{g} = \langle x, y\rangle$, $[x,y] = y$, with its representation
space $W = \langle z, t\rangle$.  Since $L$ is Frobenius $(\Delta(L) = (1+\alpha)^2 t \neq 0)$ we see that
$$i(L) = 0 = \dim W - \dim\mathfrak{g}$$
Then $L$ satisfies the Gelfand-Kirillov conjecture by Theorem 1.1 combined with Proposition 2.1 of [O6].\hfill $\square$\\
\ \\
We conclude by producing explicitly a set of Weyl generators of $D(L)$ as follows:\\
As $W = \langle z,t\rangle$ we have $R(W) = \Bbb C(z,t)$ and $R(W)^\mathfrak{g} = Z(D(L)) = \Bbb C$.  Next, we put $q_1 = z$, $q_2 = t$. By the proof of [O6, Theorem 1.1] $p_1$, $p_2$ are
the solutions of the following system of equations:
$$x = [x,q_1] p_1 + [x,q_2]p_2 \ \ \ y = [y,q_1]p_1 + [y, q_2]p_2$$
which simplifies to
$$x = \alpha zp_1 + (1 + \alpha) tp_2\ \ \ \mbox{and}\ \ \ y = tp_1$$
Hence, $p_1 = t^{-1}y$ and $p_2 = (1+\alpha)^{-1} t^{-1} (x - \alpha zt^{-1} y)$.\\
\ \\
Then, $p_1$, $p_2$, $q_1$, $q_2$ form a set of Weyl generators over $\Bbb C$ of $D(L)$ by the proof of [O6, Theorem 1.1].  Hence, $D(L) \cong D_2(\Bbb C)$.\\
\ \\
{\bf Acknowledgments}\\
\\
We are very grateful to Jacques Alev for his inspiring questions on the polynomiality of the Poisson center.  We also would like to thank Doran Shafrir and Rupert Yu for providing some essential 
information.  Special thanks go to our colleague Peter De Maesschalck for writing some efficient programs in MAPLE. Finally, we thank the referee for making some necessary corrections.\\
Part of this paper was presented at the Institut Henri Poincar\'e (Paris) and at the University of Reims.\ \\
\ \\
\renewcommand\bibname{\normalsize{References}}

\end{document}